\documentclass[ijoc]{informs3wp}

\OneAndAHalfSpacedXII 
\usepackage[pdfpagemode={UseOutlines},bookmarks=true,bookmarksopen=true,
bookmarksopenlevel=0,bookmarksnumbered=true,hypertexnames=false,
colorlinks,linkcolor={blue},citecolor={blue},urlcolor={red},
pdfstartview={FitV},unicode,breaklinks=true]{hyperref}
\hypersetup{urlcolor=blue, colorlinks=true}
\usepackage{algorithm}
\usepackage{algpseudocode}
\usepackage{tikz}

\usepackage{natbib}
 \bibpunct[, ]{(}{)}{,}{a}{}{,}%
 \def\bibfont{\small}%
\newcommand{\transpose}{{\mbox{\tiny T}}}

\newcommand{\cR}{{\mathcal{R}}}

\newcommand{\cO}{{\mathcal{O}}}

\newcommand{\bP}{\textbf{P}}

\newcommand{\bB}{\textbf{B}}

\newcommand{\bx}{\textbf{x}}

\newcommand{\bd}{\textbf{d}}

\newcommand{\bc}{\textbf{c}}

\newcommand{\bA}{\textbf{A}}

\newcommand{\ba}{\textbf{a}}

\newcommand{\bz}{\textbf{z}}
\newcommand{\bw}{\textbf{w}}

\newcommand{\br}{\textbf{r}}

\usepackage{mathtools}
\DeclarePairedDelimiter\ceil{\lceil}{\rceil}

\usepackage{mathtools}

\newif\ifnotes\notestrue
\def\tt#1{\texttt{#1}}

\def\htien#1{}

\newenvironment{proofs}[1][Proof]{\noindent\textbf{#1.} }{\hfill \ \rule{0.5em}{0.5em}}
\usepackage{bm}

\usepackage{multirow}
\usepackage{booktabs}
\usepackage{algpseudocode}
\usepackage{amsmath}
\usepackage{subcaption}

\EquationsNumberedThrough    

\TheoremsNumberedThrough     
\ECRepeatTheorems  %

\begin{document}


\RUNAUTHOR{Qian Shao, Tien Mai, and Shih-Fen Cheng}

\RUNTITLE{Constrained Pricing in Choice-based Revenue Management}

\TITLE{Constrained Pricing in Choice-based Revenue Management}

\ARTICLEAUTHORS{%
\AUTHOR{Qian Shao}
\AFF{School of Computing and Information Systems,
Singapore Management University, \EMAIL{qianshao.2020@phdcs.smu.edu.sg}}

\AUTHOR{Tien Mai}
\AFF{School of Computing and Information Systems,
Singapore Management University, \EMAIL{atmai@smu.edu.sg}}

\AUTHOR{Shih-Fen Cheng}
\AFF{School of Computing and Information Systems,
Singapore Management University, \EMAIL{sfcheng@smu.edu.sg}}

} 

\ABSTRACT{%
We consider a dynamic pricing problem in network revenue management where customer behavior is predicted by a choice model, i.e., the multinomial logit (MNL) model. The problem, even in the static setting (i.e., customer demand remains unchanged over time), is highly non-concave in prices. Existing studies mostly rely on the observation that the objective function is concave in terms of purchasing probabilities, implying that the static pricing problem with linear constraints on purchasing probabilities can be efficiently solved. However, this approach is limited in handling constraints on prices, noting that such constraints could be highly relevant in some real business considerations.  To address this limitation, in this work, we consider a general pricing problem that involves constraints on both prices and purchasing probabilities. To tackle the non-concavity challenge, we develop an approximation mechanism that allows solving the constrained static pricing problem through bisection and mixed-integer linear programming (MILP). We further extend the approximation method to the dynamic pricing context. Our approach involves a resource decomposition method to address the curse of dimensionality of the dynamic problem, as well as a MILP approach to solving sub-problems to near-optimality. Numerical results based on generated instances of various sizes indicate the superiority of our approximation approach in both static and dynamic settings.
}%





\maketitle


\section{Introduction}

In revenue management (RM) \citep{gallego1997multiproduct}, pricing is an essential strategy employed across diverse industries such as airlines, hospitality, and retail to maximize profitability. The core of this strategy involves setting optimal prices at the right times to effectively balance demand and supply, aiming to sell the right product to the right customer at the right moment. This requires sophisticated decision-making processes that take into account consumer behavior, market conditions, competitor actions, and inventory levels. Dynamic pricing strategies, which adjust prices in real-time according to market fluctuations, are particularly common, utilizing data analytics and predictive modeling to anticipate demand and modify prices, thus maximizing revenue opportunities and enhancing customer satisfaction with prices that reflect current market values.

Pricing strategies are crucial in various applications, helping to balance demand, maximize profitability, and secure competitive advantages. For instance, dynamic pricing algorithms in retail adjust to demand and supply changes to optimize revenues during peak shopping periods \citep{talluri2004theory,elmaghraby2003dynamic,Levin2009}. Airlines and hotels use revenue management systems to set prices based on booking trends and expected occupancy, maximizing earnings per available seat or room \citep{talluri2004theory,Cross1997}. In telecommunications, pricing varies with data usage, network demand, and customer loyalty to attract and retain users \citep{Audestad2007}. Energy companies employ time-based pricing to manage consumption patterns during peak and off-peak hours \citep{Borenstein2005}. Furthermore, in the digital economy, subscription models with tiered pricing structures accommodate different user needs and consumption rates \citep{Lambrecht2006}. Each application not only mirrors the overarching business strategy but also adapts to technological shifts and evolving market dynamics, underscoring the critical role of pricing in managing business operations and customer relationships.

Pricing problems, even in static settings where customer demand remains constant over time, involve solving a nonlinear optimization problem. It is well-documented that, in certain cases, such as when customer choice probabilities are linear, the pricing issue can be formulated as a convex quadratic programming problem, making it relatively straightforward to manage \citep{gallego1997multiproduct}. However, for more complex nonlinear demand functions, like the popular multinomial logit (MNL) model, the problem generally becomes non-convex in terms of pricing variables, introducing significant challenges when price constraints are taken into account. Previous research has often relied on the observation that the objective function (i.e., expected revenue as a function of the price vector) can be transformed into a concave function with respect to purchasing probabilities \citep{zhang2017dynamic,zhang2013assessing}, thus allowing the pricing problem with constraints on customer choice probabilities to be tackled through convex optimization. However, this method has limitations when there are price constraints (e.g., linear constraints that describe relationships between product prices), which are often crucial in practical scenarios. In this work, we overcome these limitations by developing a solution method that allows the (static) pricing problem to be solved to near-optimality, under constraints on both purchasing probabilities and prices.

In this paper, we explore pricing problems under MNL-based customer demand, incorporating constraints on purchasing probabilities—which can be used to model resource and inventory requirements—and price constraints that can capture additional business requirements. Examples of such price constraints include setting certain products' prices higher than others, or ensuring the total prices of some product bundles do not exceed specified lower or upper bounds. Given that the pricing problem is highly non-convex, it necessitates innovative solution methods even after typical transformations into purchasing probability space. We address these challenges by employing a piecewise linear approximation method that can achieve near-optimal solutions for the static pricing problem with arbitrarily small approximation errors. Additionally, we integrate our approach with a dynamic programming decomposition technique to tackle dynamic pricing problems. Our contributions are detailed as follows:
\begin{enumerate}
    \item \textbf{Static Pricing with Constraints:} We begin by tackling a static pricing problem that involves constraints on both prices and purchasing probabilities. To address the inherent non-convexity of this problem, we develop a solution method that combines bisection and piecewise linear approximation. Specifically, we employ Dinkelbach's transform to convert the optimization problem into a sequence of sub-problems. Each sub-problem is then approximated by a mixed-integer linear program (MILP), solvable by commercial solvers like CPLEX or GUROBI. By analyzing the convexity and concavity within these sub-problems, we identify opportunities to safely relax certain binary variables within the MILP approximation, significantly enhancing its performance. This strategy enables us to solve the static pricing problem to any desired precision.
    
    \item \textbf{Analysis of Approximation Errors:} We rigorously examine the approximation errors introduced by our solution method. Specifically, we demonstrate that solving the approximated problem (through a series of MILPs) results in objective values that are within an \(\mathcal{O}(1/K)\) neighborhood of the optimal values, where $K$ is the number of breakpoints used to discretize the prices.
    
    \item \textbf{Dynamic Pricing Adaptation:} We extend our discretization approach to solve dynamic pricing, employing the dynamic programming decomposition approach to handle the curse of dimensionality. Our approach is heuristic in this context, but it offers the advantage of providing near-optimal solutions for each backward step of the dynamic programming decomposition.
    
    \item \textbf{Numerical Analysis:} We conduct extensive comparisons with standard baselines in the literature, including a standard nonlinear solver, an adaptation of the convex transformation commonly used in previous studies, and Gurobi's piecewise linear approximation methods that incorporate state-of-the-art techniques for nonlinear optimization. Our results clearly demonstrate the superiority of our approaches for both static and dynamic pricing problems, highlighting significant improvements over existing methods.
\end{enumerate}

\noindent
\textbf{Paper Outline:}
The paper is organized as follows: Section \ref{sec:review} provides a literature review. Section \ref{sec:formulation} introduces the problem formulations for both static and dynamic pricing. Section \ref{sec:SP} discusses how to solve the constrained static pricing problem, while Section \ref{sec:DP} shows how these techniques are applied to the constrained dynamic pricing problem. Section \ref{sec:Numerical Results} provides numerical comparisons and analyses. Finally, Section \ref{sec:concl} concludes the paper. The appendix includes proofs not contained in the main text and some additional  details.

\noindent
	\textbf{Notation:}
	Boldface characters represent matrices or vectors or sets, and $a_i$ denotes the $i$-th element of $\ba$ if $\ba$ is indexable. We use $[m]$, for any $m\in \mathbb{N}$, to denote the set $\{1,\ldots,m\}$.

\section{Literature Review}\label{sec:review}
Our work is closely relevant to the literature on pricing in RM, approximate dynamic programming (ADP), and piecewise linear approximation. The comprehensive overview of the RM literature can be found in \cite{talluri2004theory}.

In the context of static pricing, the study conducted by \cite{hanson1996optimizing} demonstrates that the profit function of multiple differentiated substitutable products, when analyzed under the MNL model, does not exhibit joint \textit{concavity} with respect to the price vector. Consequently, direct application of concave optimization techniques cannot yield the optimal solutions. A similar observation is made in the research conducted by \cite{dong2009dynamic}. Therefore, instead of utilizing the price vector as decision variables for static pricing, an alternative approach involves employing the purchasing probabilities as decision variables. relying a one-to-one correspondence with the price vector  (see \cite{xue2007demand}, \cite{dong2009dynamic} and \cite{zhang2013assessing} for the MNL model, and \cite{li2011pricing} for the Nested Logit (NL) model). In these works, the authors successfully utilize the transformation into the purchasing probability space, effectively converting the pricing problem into a convex one. However, when pricing constraints are included  (such as simple linear constraints on prices), the previously mentioned transformation becomes non-convex, leading to potentially sub-optimal solutions. Therefore, our approach represents a significant advancement in efficiently solving static pricing problems with complex business constraints.

Dynamic pricing, widely recognized as a crucial subarea of revenue management, has witnessed substantial and remarkable growth in recent years \citep{zhang2013assessing}. Dynamic pricing involves the strategic determination of optimal selling prices for products or services in a flexible pricing environment. This concept applies to various scenarios, including online vendors selling through the Internet and physical stores utilizing digital price tags. Dynamic pricing techniques have become pervasive in modern business practices and are often regarded as an essential component of pricing strategies. There are several comprehensive and extensive reviews of the extensive literature on dynamic pricing, including the  review articles: \cite{bitran2003overview}, \cite{elmaghraby2003dynamic}, \cite{deksnyte2012dynamic}, \cite{gonsch2009dynamic}.

Dynamic pricing presents a formidable challenge due to its huge state space, and ADP approaches have become popular for tackling this challenge. Within the realm of ADP for RM, \cite{de2003linear} explore methodologies based on linear programming, while \cite{adelman2007dynamic} pioneered the use of linear programming for approximate dynamic programming in RM. This involves approximating the value function with a linear combination of basis functions and integrating this approximation into the linear programming formulation of the dynamic equations. \cite{zhang2009approximate} utilize a linear functional approximation to generate dynamic bid prices for use in dynamic programming decomposition. \cite{kunnumkal2008refined} introduce a refined deterministic linear programming model for network RM, applying a Lagrangian relaxation to relax capacity availability constraints. \cite{kunnumkal2010new} develop an alternative dynamic programming decomposition method to overcome some limitations of the method used in \cite{liu2008choice}, selecting revenue allocations through an auxiliary optimization that accounts for the probabilistic nature of customer choices. \cite{topaloglu2009using} applies a Lagrangian relaxation approach to a dynamic programming formulation of the network RM problem. More recently, \cite{erdelyi2011using} adopt similar techniques for dynamic pricing issues within a network context. \cite{meissner2012network} propose an approximate programming approach for network RM models incorporating customer choice, approximating the value function of the Markov decision process using a non-linear function that exhibits resource inventory separability. \cite{zhang2013assessing} show that dynamic programming decomposition offers a revenue upper bound tighter than that from deterministic methods. \cite{zhang2011improved} explore a nonlinear, nonseparable functional approximation for representing the value function in a dynamic programming formulation of network RM with customer choice, formulating it as a nonlinear optimization problem with complex constraints. Lastly, \cite{zhang2017dynamic} conduct a comprehensive numerical study comparing several heuristic approaches from the literature, including deterministic linear programming with resolving and three variants of dynamic programming decomposition. To the best of our knowledge, all the aforementioned studies and other relevant papers that employ the MNL customer demand model typically overlook the constraints on prices, focusing on the use of the purchasing-probability transformation to solve each step of the dynamic programming decomposition. In contrast, our work introduces an advanced approach to handle price constraints in dynamic pricing. Our subsequent numerical experiments demonstrate that our approach provides significantly better pricing solutions when solving large-scale dynamic pricing problems.

Our proposed methods utilize a piece-wise linear approximation approach to transform the nonlinear program into a sequence of MILPs. The literature on piece-wise linear approximation is extensive \citep{PWLA_lin2013review,PWLA_lundell2013refinement,PWLA_westerlund1998extended,PWLA_lundell2009some}, and such techniques have been incorporated into some advanced solvers like GUROBI for solving mixed-integer nonlinear programs \citep{GUROBI}. To the best of our knowledge, existing studies and solvers primarily focus on univariate functions of specific types, such as exponential, sigmoidal, or quadratic functions, which are not directly applicable to the fractional structure of the pricing problem. In our experiments, we also provide a comparison that demonstrates our approximation method outperforming an approach that employs GUROBI's piece-wise linear approximation  to solve each step of the bisection.



\section{Problem Formulation}\label{sec:formulation}
We present mathematical formulations for both static and dynamic pricing problems where customer choice behavior is predicted by a MNL model.

\subsection{Static Pricing}

Let us consider  a pricing problem with $n$ products, indexed as $\{1,2,...,n\}$.  The non-purchasing option is considered a special product and indexed as ``0''.  Let $r_j$ be the price (i.e., revenue) of  product $j \in [n]$  and $P_j(\br)$ be the probability of that product $j\in [n] \cup \{0\}$   is chosen. Under resources constraints, it is assumed that there is a network with  $m$ resources, where each resource $i\in [m]$ has a capacity $c_i>0$. The resources can be combined to produce $n$ products.  An $m\times n$ matrix $\bA$ is used to represent the resource consumption, where the $(i,j)$-th element $a_{ij}$ denotes the quantity of resource $i$ consumed by one unit of product $j$; $a_{ij} =1$ if resource $i$ is used by product $j$ and $a_{ij}=0$ otherwise. Let $\bA_i$ be the $i$-th row of $\bA$ and $\bA^j$ be the $j$-th column of $\bA$, respectively. To simplify the notation, we use $j \in \bA_i$ to indicate that product $j$ uses resource $i$ and $i\in \bA_j$ to indicate that resource $i$ is used by product $j$.
It is then required that the total consumption of each resource cannot exceed the resource capacity. The resource  constraints can be stated as $\bA \bP(\br) \leq \bc$, where $\bP(\br)$ is the vector of all the purchasing probabilities $\bP(\br) = (P_j(\br);~ j\in [n])$. 

The MNL model is one of the most popular choice models widely used to model and predict people demand \citep{train2009discrete} and  has also been intensively used in the RM literature \citep{Talluri2004,dong2009dynamic} . This model is based on the random utility maximization (RUM) framework, which assume that each product $j\in [n]$ is associated with a deterministic utility $v_j$. To capture the effect of the prices $\br$ on the choice probabilities, these utilities  are typically modelled as a linear function of the prices, i.e., $v_j = (a_j-r_j)/b_j$, $j\in [n]$, where $a_j$ and $b_j$ are some choice parameters that can be inferred from fitting the choice probabilities with  observations of how customers made choices. In particular,  $1/b_j$ is typically referred to as a price sensitivity parameter  associated with product $j\in [n]$ and should take a positive value, which is because  higher prices should yield lower utility values. Under the MNL choice model, the probability of purchasing product $j\in [n]\cup \{0\}$ can be computed as:
\begin{align}
\label{purchasing probability}
P_j(\br)= \frac{e^{v_j}}{\sum_{j'\in [n]\cup \{0\}} e^{v_{j'}}}  = \frac{e^{(a_j - r_j)/b_j}}{\sum_{j\in [n]} e^{(a_j - r_j)/b_j} + e^{v_0}}
\end{align}
where $v_0$ represents the utility of  the non-purchasing option. Without the loss of generality and for ease of notation, let us assume that $v_0 = 0$. Under the MNL model, the objective function $f(\br)$ is no quasi-concave  (Hanson and Martin (1996)). The constraint $\bA\bP(\br)\leq \bc$ is also not convex, making the pricing problem non-convex within the price space and challenging to solve. Nevertheless, it can be shown that there is one-to-one mapping between a  price vector $\br$  and a purchasing probability vector $\bP$, which allows to  write the prices as a well-defined function of the purchasing probabilities, which we can denote as $\br(\bP)$. Consequently,  the objective function, in the $\bP$ space,  can be written as $\br(\bP)^\transpose \bP$. Under the assumption that the price-sensitivity parameters $1/b_j$ are the same over all the products, this objective function is concave  in $\bP$, making the pricing  problem  convex within the $\bP$ space. A majority of prior works  relies on this observation to solve the pricing problem. It is worth noting that this results hold for any choice model in the General Extreme Value (GEV) facility of  choice models \citep{zhang2018multiproduct}, which includes MNL and other popular choice models such as the nested logit and cross-nested logit models.  

It's noteworthy that assumption that the price-sensitivity parameters are the same  across product, i.e., $b_j=b_{j'}, ~\forall j,j'\in [n]$, is  a crucial condition for transforming the pricing problem into a convex  optimization problem \citep{dong2009dynamic,zhang2018multiproduct}. However, this assumption would be restricted in numerous real-world scenarios. Moreover, the convexity does not hold if  the problem involves constraints on $\br$. In these regards, our framework presents several advantages. Specifically, it eliminates the need for homogeneity of $b_j$ across products $j\in [n]$, allowing the price-sensitivity parameters to take arbitrary positive values. Additionally, it accommodates the incorporation of any linear constraints on the prices $\br$.  Particularly, our work addresses the following pricing problem with resource  and price constraints:  
\begin{align}
 \textbf{[SP]} ~~~~~ \max_{\br} &\left\{ F(\br) =  \sum_{j=1}^n r_j P_j(\br)  \right\} \label{prob:SP-2} \\
    \mbox{subject to} &{\quad \bA \bP(\br) \leq \bc  }\nonumber \\
    &{\quad \bB \br \leq \bd  }\nonumber \\
    &\quad l_j\leq r_j \leq u_j,~~\forall j\in [n],\nonumber
\end{align}
where $l_j$ and $u_j$ are lower  and upper bounds of the price decision $r_j$, and $\bB\br\leq \bd$ are other constraints on $\br$, reflecting other business requirements on the price decisions. Relevant practical constraints may include the following:
\begin{itemize}
\item Certain item prices must be reduced compared to previous pricing decisions, often due to considerations such as  promotions. Such constraints can be formulated as $r_j < \overline{r}_j$ for a given $j$, where $\overline{r}_j$ represents the previous price of item $j$.
\item Some items needs to have higher prices (or lower) than some others: $r_j \geq r_{j'}$ for some $j,j'\in [n]$.
\item The cumulative price of a bundle of items should not surpass a specified threshold: $\sum_{j\in [n]} \alpha_j r_j\leq\beta$, where $\alpha_j\geq 0$ and $\beta >0$ are predefined parameters.
\item  The new prices should not deviate significantly from a previous set of prices, which can be expressed as: $\sum_{j\in [n]} |r_j - \overline{r}_j| \leq \beta$, where $\overline{r}_j$, for $j \in [n]$, represents prior pricing values.
{Incorporating the aforementioned constraints would undoubtedly enhance the practical applicability of the pricing problem formulations and optimization methods in real-world business scenarios.}
\end{itemize}

\subsection{Dynamic Pricing}

In  dynamic pricing,  a network with $m$ resources is distributed over products in a selling horizon. We divide the selling horizon into $T$ time periods, and indexed by $t=1,...,T,T+1$, where $t=1$ and $t=T+1$ represent the beginning and the end of the selling horizon. In each time period $t$, the probability of one customer arrival is $\lambda$, and the probability of no customer arrival is $1-\lambda$. We assume that a customer arriving at time $t$ only buy one product. Also, let  $\bx$ be the vector of remaining resource capacity at time $t$, $v_t(\bx)$ be the maximum expected revenue given state $\bx$ at time $t$, $\br_t$ is the pricing decision at time $t$, and $R_t(\bx)$ is  the feasible set of the pricing decision $\br_t$ at time $t$.   The dynamic pricing problem can be formulated as follows:
\begin{equation}
    \textbf{[DP]} ~~~~\max_{\substack{\br_t}} \left\{v_1(\textbf{c})\right\}
\end{equation}
where the value function $v_t(\bx)$, $t = 1,\ldots,T+1$, satisfy the  following  Bellman equation:
\begin{align}
     v_t(\bx) & = \max \limits_{\br_t \in R_t(\bx)} \left\{ \sum_{j=1}^n \lambda P_j(\br_t)[r_{t,j} + v_{t+1}(\bx-\bA^{j})] + (\lambda P_0(\br_t) +1-\lambda)v_{t+1}(\bx)\right\}  \label{Bellman Equation} \\
    &= \max \limits_{\br_t \in R_t(\bx)} \left\{ \sum_{j=1}^{n} \lambda P_j(\br_t)[r_{t,j} + v_{t+1}(\bx-\bA^{j}) - v_{t+1}(\bx)] \right\} + v_{t+1}(\bx) \nonumber\\
    &= \max \limits_{\br_t \in R_t(\bx)} \left \{\sum_{j=1}^{n}\lambda P_j(\br_t)[r_{t,j} - \triangle_{j}v_{t+1}(\bx)] \right \} + v_{t+1}(\bx) \nonumber 
\end{align}
Where $\triangle_{j}v_{t+1}(\bx)=v_{t+1}(\bx)-v_{t+1}(\bx-\bA^{j})$ represents the expected revenue of selling one product $j \in [n]$ at time $t+1$. The boundary conditions are  $v_{T+1}(\bx) = 0,\forall x $ and $v_t(\textbf{0}) = 0, \forall t\in [T] $. In the above, $R_t(\bx) \subset \Re^n_{+}$ , and  $R_{t}(\bx) = \{ \br^\infty \}$ if  $\bx\geq \bA^j$. The price $\br^\infty$ is known as the null price \citep{gallego1997multiproduct}, characterized by the property that $P_j(\br)=0$ when $r_j= r^\infty_j$. Consequently, in situations where there is an insufficient quantity of resources to meet the demand for product $j \in [n]$, the demand can effectively be eliminated by setting $r_j=r_j^\infty$.

At each time step $t\in [T]$, the feasible set $R_t(\bx)$ can be represented by the following set of constraints.
\[
R_t(\bx) = \left\{\br \in \Re_+^n\Bigg\vert\,
  \begin{array}{@{}cc@{}}
    l_j\leq r_j\leq u_j,~~\forall j\in [n]
      &\quad (c1) \\
    \lambda \bA \bP(\br) \leq \bx &\quad (c2) \\
    \bB \br \leq \bd &\quad (c3)
  \end{array}
  \right\}
\]
where $(c1)$ are box constraints  on the prices $\br$,  $(c2)$ is to ensure that the expected resource consumption does not exceed the  current resource capacity $\bx$, and constraints $(c3)$ are some  linear constraints  on the prices. Here we note that the constraints in $(c1)$ and $(c3)$ can  be time-dependent but we keep them independent of time $t$ for notational simplicity.

The dynamic problem [\textbf{DP}] can  be solve by  iteratively computing $v_t(\bx)$ for  $t = T ,... 1$.  Each iteration requires to solve a static pricing problem of the form:
\[
v_t(\bx) = \max \limits_{\br \in R_t(\bx)} \left \{\sum_{j \in  [n]} \lambda P_j(\br)\left[r_{j} - \triangle_{j}v_{t+1}(\bx)\right] \right \} + v_{t+1}(\bx) \nonumber 
\]
Even though this classical approach enables the exact computation of $v_t(\bx)$, it becomes impractical in high-dimensional settings. 
In practice, it is a common approach to decompose the network problem into a collection of smaller sub-problems, with each sub-problem focusing on a single resource, such as a single leg or a single stay night. Such an approach can be found in \cite{talluri2004theory},~\cite{zhang2013assessing},~\cite{zhang2011improved},~\cite{zhang2017dynamic}. While we employ a similar  decomposition method to tackle our dynamic pricing  problem, the inclusion of constraints on the price decisions make prior methods no-longer suitable. Specifically, in previous works, the convexity of the static pricing problem within the $\bP$ space has been leveraged to solve sub-problems. However, in our context, where convexity no longer holds, we introduce a novel approach using bisection (to be discussed in the following sections) and mixed-integer linear programming. This enables us to address sub-problems more efficiently. Additionally, as highlighted in the preceding section, our approach is not bound by the price-homogeneity assumption, enhancing the applicability of the pricing formulation in broader real-world scenarios.

\section{{Constrained Static Pricing}}\label{sec:SP}
\label{section SP}

In this section, we discuss our solution method for the constrained pricing problem in \textbf{[SP]}. To tackle the highly nonlinear problem, we will leverage  the Dinkelbach transform to simplify the fractional structure, and a piece-wise linear approximation approach to approximate the nonlinear problem by a MILP, for which an off-the-shelf solver can be conveniently used.

\subsection{Piece-wise Linear Approximation}

We first present our approach to approximate the nonlinear components of the objective function of \textbf{[SP]} by  piece-wise linear functions, which will later allow us to solve the nonlinear pricing problem by bisection and MILP solvers. First, for ease of notation, let us denote:
\begin{align*}
    f_j(r_j) = r_je^{(a_j-r_j)/b_j};~~
    g_j(r_j) = e^{(a_j-r_j)/b_j}    
\end{align*}
To approach each nonlinear univariate function $f_j(r_j)$, $g_j(r_j)$ by a piece-wise linear function, we discretize the interval $ [l_j,u_j]$ into $K$ sub-intervals of equal sizes and  approximate these functions as
\begin{align*}
    f_j(r_j) \approx \widehat{f}_j(r_j) = f\left(k \Delta_j\right) + 
\gamma^f_{jk}\left(r_j - k\Delta_j\right);~~
 g_j(r_j) \approx \widehat{g}_j(r_j) = g\left(k \Delta_j\right) + 
\gamma^g_{jk}\left(r_j - k\Delta_j\right)
\end{align*}
where $k = \left\lfloor \frac{r_j}{\Delta_j}\right\rfloor$ and $\Delta_j = \frac{u_j-l_j}{K}$ and $\gamma^f_{jk}$ is the slope of the linear segments between two discretization points $((k\Delta_j, f(k\Delta_j)); ((k+1)\Delta_j, f((k+1)\Delta_j)))$, and  $\gamma^g_{jk}$ is the slope of the linear segments between two discretization points $((k\Delta_j, g(k\Delta_j)); ((k+1)\Delta_j, g((k+1)\Delta_j)))$, computed as 
\begin{align*}
    \gamma^f_{jk} = \frac{f((k+1)\Delta_j) - f(k\Delta_j)}{\Delta_j};~ \gamma^g_{jk} &= \frac{g((k+1)\Delta_j) - g(k\Delta_j)}{\Delta_j}
\end{align*}
We now consider the following approximate problem, obtained by replacing $f_j(r_j)$ and $g_j(r_j)$ with $\widehat{f}_j(r_j)$, $\widehat{g}_j(r_j)$, respectively, in \textbf{[SP]}
\begin{align}
 \textbf{[ASP]} ~~~~~ \max_{\br} &\left\{ \widehat{F}(\br) =  \frac{\sum_{j\in [n]} \widehat{f}_j(r_j)}{ 1+ \sum_{j\in [n]} \widehat{g}_j(r_j)} \right\} \label{prob:Approx-SP-1} \\
    \mbox{subject to} &{\quad \bA \widehat{\bP}(\br) \leq \bc  }\nonumber \\
    &{\quad \bB \br \leq \bd  }\nonumber \\
    &\quad l_j\leq r_j \leq u_j,~~\forall j\in [n].\nonumber
\end{align}
where $\widehat{\bP}(\br)$ is a vector of size $n$ with entries 
\[
\widehat{P}_j(\br) =\frac{\widehat{g}_j(r_j)}{1 + \sum_{j\in [n]} \widehat{g}_j(r_j) },~\forall j\in [n].
\]
Having this, the constraints $\bA\widehat{\bP}(\br)\leq \bc$ can be rewritten equivalently as:
\[
\sum_{j\in [n]} (a_{ij}-c_i)\widehat{g}_j(r_j) \leq c_i,~~~\forall i\in [m].
\]

\subsection{Solving the Approximate Problem}

We discuss how to efficiently solve the approximate problem \textbf{[ASP]}. Our approach involve converting \textbf{[ASP]} into to  bisection procedure  where at each step requires to solve a subproblem, which is  a mixed-integer nonlinear program. We then show that the subproblem can be reformulated as a MILP. To facilitate our later exposition, let us first denote by $\widehat{\cR}$ be the feasible set of the approximate problem \textbf{[ASP]}
\[
\widehat{\cR} = \left\{\br \in \Re_+^n\left\vert\,
  \begin{array}{@{}ll@{}}
    l_j\leq r_j\leq u_j,~~\forall j\in [n]
      & \\
    \sum_{j\in [n]} (a_{ij}-c_i)\widehat{g}_j(r_j) \leq c_i,~~~\forall i\in [m] & \\
    \bB \br \leq \bd &
  \end{array}
  \right.\right\}
\]

\noindent \textbf{{Bisection.}}  We use the Dinkelbach transform \citep{dinkelbach1967nonlinear} to convert the fractional objective function of \textbf{[ASP]} into a non-fractional one.  The Dinkelbach transform works by observing  that \textbf{[ASP]} can be equivalently converted into the  following univariate optimization problem:
\begin{align*}
     \max \left\{\delta\Big|~ \exists \br \in \widehat{\cR}\texttt{ s.t. } \widehat{F}(\br) \geq \delta\right\}
\end{align*}
In other words, we can solve \textbf{[ASP]}  by finding the highest value of a threshold $\delta$ such that there exists some feasible price sets $\br$ such that $\widehat{F}(\br) \geq \delta$, which is equivalent to:
\[
\sum_{j\in [n]} \widehat{f}_j(r_j) \geq \delta \left(1+\sum_{j\in[n]} \widehat{g}_j(r_j)\right).
\]
The highest possible value of $\delta$ can be computed by a binary search where in each round of the search the feasibility problem stated above is solved for a particular $\delta$. For given $\delta$, the feasibility problem is easily solved by checking if the maximum of following optimization is non-negative or not
\begin{equation}
\label{obj_st}
 \textbf{[BOPT]} \quad \quad  \max\limits_{\br\in\widehat{\cR}}~ \left\{ G(\br,\delta) = \sum_{j\in [n]} \widehat{f}_j(r_j) - \delta \left(1+\sum_{j\in[n]} \widehat{g}_j(r_j)\right) \right\}
\end{equation}
The binary search procedure  over $\delta$ is shown in Algorithm \ref{Algorithm binary search}. The lower and upper bounds are chosen in such a way that the range $[L,U]$ should cover the optimal value of \textbf{[ASP]}. We can simply choose $L=0$ and  $U=\max_{j\in [n]} u_j$, due to the fact that
\[
F(\br) \leq \max_{j\in [n]} \{ r_j\} \leq  \max_{j\in [n]} \{u_j\}.
\]
It can be shown that, if \textbf{[BOPT]}  is solved to optimality at each step, then the binary search will always terminate after  about $\log\left(\frac{U-L}{\xi}\right)$  and return  a solution within a $\xi$ neighbourhood of optimality. We state this result in Proposition [\ref{prop:bisection-gap}] below.
\begin{proposition}\label{prop:bisection-gap}
    For any $\xi\geq 0$, Algorithm \ref{Algorithm binary search} always terminates after no more than $\ceil{\log\left(\frac{U-L}{\xi}\right)}$ steps, and return a solution $\widehat{\br}$ such that $\widehat{F}(\widehat{\br}) \geq \widehat{F}(\br^*) -\xi$, where $\br^*$ is an optimal solution to \textbf{[ASP]}. 
\end{proposition}



\begin{algorithm}
\caption{Binary Search Template}
\label{Algorithm binary search}
\begin{algorithmic}[1]
\State $U = \max \limits_{j\in [n]} u_j$, $L = 0$, $\xi > 0$
\While{$U - L \ge \xi$}
    \State $\delta = (U + L)/2$
    \State $G^* = \max_{\mathbf{r} \in \widehat{\mathcal{R}}} \{G(\mathbf{r}, \delta)\}$
    \State \text{If} $G^* \ge 0$ \text{then} $L = \delta$ \text{, else} $U = \delta$
    
\EndWhile
\State \textbf{return} $\delta$, $\mathbf{r}$ from the last solution of \textbf{[BOPT]}.
\end{algorithmic}
\end{algorithm}

In the following we show that \textbf{[BOPT]} can be converted into a MILP. 

\noindent \textbf{{MILP Formulation for \textbf{[BOPT]}}}. We show that the subproblem of the binary search procedure discussed above can be reformulated as a MILP. To this end, let us introduce binary variable $z_{jk}$    and continuous variable $w_{jk}$, $j\in [n]$, $k\in \{0,\ldots,K-1\}$ such that 
 \begin{align*}
     &z_{jk} \geq z_{j,k+1}, ~ k =0,\ldots, K-2\\
     &z_{jk}\leq w_{jk} \leq 1,~ k =0,\ldots, K-1\\
     &w_{j,k+1} \leq z_{jk}, ~ k =0,\ldots, K-2 
 \end{align*}
Here, the binary variables $z_{jk}$ are to capture the part $\left\lfloor \frac{r_j}{\Delta_j}\right\rfloor \Delta_j$ of the price $r_j$  and the continuous variables $w_{jk}$ are to represent the remaining part $r_j - \left\lfloor \frac{r_j}{\Delta_j}\right\rfloor \Delta_j$. It can be seen that there is a unique mapping between any price value $r_j$ with a set $\{z_{jk}, w_{jk}|~ k=0,\ldots K-1\}$ that satisfies the above conditions such that $r_j = \Delta_j \sum_{k = 0}^{K-1}w_{jk}$. Then, the approximate functions $\widehat{f}_j(r_j)$ and $\widehat{g}_j(r_j)$ can computed via $z_{ik}$ and $w_{jk}$ as
\begin{align*}
    \widehat{f}_j(r_j) = f_j(l_j) + \Delta_j\sum_{k=0}^{K-1}\gamma^f_{jk} w_{jk};~~
    \widehat{g}_j(r_j) = g_j(l_j) + \Delta_j\sum_{k=0}^{K-1}\gamma^g_{jk} w_{jk}
\end{align*}
As a result,  \textbf{[BOPT]} can be formulated equivalently as the following MILP:
\begin{align}
\textbf{[MILP]} ~~~~~ \max\limits_{\bw,\bz,\br} & \sum_{j\in [n]} (f_j(l_j)-\delta g_j(l_j)) + \sum_{j \in [n]} \sum_{k = 0}^{K-1}\Delta_j(\gamma^f_{jk}-\delta \gamma^g_{jk})w_{jk}-\delta  \label{prob:MILP} \tag{\sf MILP}\\
\mbox{subject to} \quad &   z_{jk}\geq z_{j,k+1}  \quad k= 0,\ldots,K-2, j \in [n] \label{ctr-milp-1} \\
& {z_{jk}} \leq w_{jk} \leq 1, \quad k = 0,\ldots,K-1, j \in [n]  \label{ctr-milp-2} \\
& w_{j,k+1} \leq z_{jk}  \quad k = 0,\ldots,K-2, j \in [n]  \label{ctr-milp-3} \\
&r_j = \Delta_j\sum_{k=0}^{K-1} w_{jk},\quad \forall j\in [n]\label{ctr-milp-4}\\
& \sum_{j\in [n]} (a_{ij} - c_i) \left(g_j(l_j) + \sum_{k=0}^{K-1} \gamma^g_{jk} w_{jk}\right) \leq c_i \quad \forall i \in [m] \label{ctr-milp-5} \\
&\bB\br \leq \bd \label{ctr-milp-6}\\
& z_{jk} \in \{0,1\}, w_{jk} \in \mathbb{R} \quad \forall k\in[K], j \in [n] \label{ctr-milp-7}
\end{align}

\noindent \textbf{Improved MILP Formulation for \textbf{[BOPT]}.} The MILP formulation in \eqref{prob:MILP} involves  $K\times n$ additional binary variables, which would be a concern when $n$ and $K$ are large. In the following we show that, be leveraging the convexity of both $g_j(\cdot)$  and $f_j(\cdot)$, one can significantly reduce the number  of binary variables in \eqref{prob:MILP}, leading to an improved MILP formulation. Let us start by introducing the following two lemmas.
\begin{lemma}\label{lm:1}
          For any $j\in [n]$,  $g_j(r_j)$ is strictly convex in $r_j$, and $f_j(r_j)$ is concave for $r_j\leq 2b_j$, and convex when $r_j\geq 2b_j$
\end{lemma}

The following lemma further shows a relation between the consecutive slope values $\gamma^f_{jk}$ and $\gamma^g_{jk}$  as $k$ increases.

\begin{lemma}\label{lm:2}
    For any $j\in [n]$, $\gamma^g_{jk}< \gamma^g_{j,k+1}$ for all $k=0,\ldots,K-1$. Moreover,   $\gamma^f_{jk}\leq \gamma^f_{j,k+1}$ when $k+2\leq \lfloor\frac{2b_j}{\Delta_j}\rfloor $ and $\gamma^f_{jk}\geq \gamma^f_{j,k+1}$ for $k\geq\lceil\frac{2b_j}{\Delta_j}\rceil $
\end{lemma}

Leveraging the properties stated in Lemma \ref{lm:2}, the following theorem shows that, without the resource constraints, part of the additional binary variables in \eqref{prob:MILP} can be relaxed. Especially, if $u_j\leq 2b_j$ for all $j\in [n]$, then \eqref{prob:MILP} can be safely converted to a linear program.

\begin{theorem}\label{th:relax-MILP}
In \eqref{prob:MILP} without the resource constraints,   binary variables  $z_{jk}$ can be relaxed for all $j,k$ such that $k\leq  \lfloor\frac{2b_j}{\Delta_j}\rfloor - 1$. That is  \eqref{prob:MILP} is equivalent to the following relaxed program:
\begin{align}
  \mathop{max}\limits_{\bw,\bz,\br} &~~~~~~~~ \sum_{j\in [n]} (f_j(l_j)-\delta g_j(l_j)) + \sum_{j \in [n]} \sum_{k = 0}^{K-1}\Delta_j(\gamma^f_{jk}-\delta \gamma^g_{jk})w_{jk}-\delta  \label{prob:MILP-2} \tag{\sf Re-MILP}\\
\mbox{subject to} \quad & ~~~~~~~~ \text{Constraints}~~ \eqref{ctr-milp-1}-\eqref{ctr-milp-4} ~~\&~~ \eqref{ctr-milp-6}\\
& ~~~~~~~~ \bm{z_{jk} \in [0,1]},~ 0\leq k\leq   \left\lfloor\frac{2b_j}{\Delta_j}\right\rfloor - 1, \text{ and } z_{jk} \in \{0,1\},~ k\geq    \left\lfloor\frac{2b_j}{\Delta_j}\right\rfloor \\
& ~~~~~~~~ \bw \in [0,1]^{n\times K} \nonumber
\end{align}
As a result, if  $u_j \leq 2b_j$ for all $j\in [n]$, then  \textbf{[BOPT]} is equivalent to the following  linear program:
     \begin{align}
  \mathop{max}\limits_{\bw,\bz,\br} &~~~~~~~~ \sum_{j\in [n]} (f_j(l_j)-\delta g_j(l_j)) + \sum_{j \in [n]} \sum_{k = 0}^{K-1}\Delta_j(\gamma^f_{jk}-\delta \gamma^g_{jk})w_{jk}-\delta  \label{prob:LP} \tag{\sf LP}\\
\mbox{subject to} \quad & ~~~~~~~~ \text{Constraints}~~\eqref{ctr-milp-1}- \eqref{ctr-milp-4} ~~\&~~ \eqref{ctr-milp-6}\\
& ~~~~~~~~ \bz,\bw\in [0,1]^{n\times K} \nonumber 
\end{align}
     
\end{theorem}
The proof can be found in the appendix. Intuitively, we leverage the result that the binary variables \( z_{jk} \) can be relaxed if \( f_j(r_j) - \delta g_j(r_j) \) is concave in \( r_j \). Theorem \ref{th:relax-MILP} is then established by identifying intervals where \( f_j(r_j) - \delta g_j(r_j) \) is concave in \( r_j \). In particular, when such intervals cover the lower and upper bounds of \( r_j \), all the additional binary variables \( z_{jk} \) can be safely relaxed and [\textbf{BOPT}] can be solved via linear programming. 
   
   \subsection{Approximation Errors}

In this section, we explore the approximation errors resulting from solving the approximated problem \textbf{[ASP]}. Specifically, we aim to answer the questions regarding how solutions to \textbf{[ASP]} approach optimal solutions to the original pricing problem \textbf{[SP]} as the number of discretization points $K$ increases. 
Firsts, let us introduce the following definition:
\begin{definition}
\textit{     A pricing solution $\br$ is said to be $\epsilon$-feasible to \textbf{[SP]}  if (and only if) 
     \[
\br \in {\cR}^{\epsilon} \stackrel{def}{=} \left\{\br \in \Re_+^n\left\vert\,
  \begin{array}{@{}ll@{}}
    l_j\leq r_j\leq u_j,~~\forall j\in [n]
      & \\
    \sum_{j\in [n]} (a_{ij}-c_i){g}_j(r_j) \leq c_i +\epsilon,~~~\forall i\in [m] & \\
    \bB \br \leq \bd &
  \end{array}
  \right.\right\}
\]
 }
   
\end{definition}
The following lemma provides an upper bound for the gap between the approximate function \( \widehat{F}(\br) \) and the original function \( F(\br) \) for all \( \br \) in its feasible set.
\begin{lemma}\label{lm:bound-F-hatF}
    For  any $\br \in \cR$, we have 
     $ \left|F(\br) - \widehat{F}(\br)\right| \leq \frac{\omega}{K},$
  where 
  \begin{align}
    \alpha &= \sum_{j\in [n]}  e^{(a_j-l_j)/b_j}\left(1+\frac{u_j}{b_j}\right)(u_i-l_i);~~
    \beta= \sum_{j\in [n]}  \frac{e^{(a_j-l_j)/b_j}}{b_j}(u_i-l_i)\nonumber\\
    \omega &= \frac{\alpha +\beta r^u}{1+\sum_{j\in [n]} e^{(a_j-u_j)/b_j}} \nonumber
  \end{align}
\end{lemma}

We are now ready to bound the approximation errors resulting from solving the approximate problem \textbf{[ASP]}. Let us first denote
\begin{align*}
\psi_i(\br) &=  \sum_{j\in [n]} (a_{ij}-c_i){g}_j(r_j), ~~
\widehat{\psi}_i(\br) =   \sum_{j\in [n]} (a_{ij}-c_i)\widehat{g}_j(r_j), ~~~\forall i\in [m]
\end{align*}
and introduce the following bound for the gap between $\psi(\br)$ and $\widehat{\psi}(\br)$
\begin{lemma}\label{lm:bound-psi}
    For any $\br  \in [l_j,u_j]^n$, we have
    $|\psi_i(\br) - \widehat{\psi}_i(\br)| \leq \frac{\eta}{K},$
    where $\eta = \sum_{j\in [n]}|a_{ij}-c_i| \frac{e^{(a_j-l_j)/b_j}}{b_j}({u_j-l_j}).$
\end{lemma}
We are now ready to bound the approximation errors resulting from solving [\textbf{ASP}]. Since establishing direct bounds for \textbf{[ASP]} is challenging, we propose a slight adaptation of \textbf{[ASP]} as follows:
\begin{align}
~~~~~ \max_{\br} &\left\{ \widehat{F}(\br) =  \frac{\sum_{j\in [n]} \widehat{f}_j(r_j)}{ 1+ \sum_{j\in [n]} \widehat{g}_j(r_j)} \right\} \label{prob:Approx-SP-2}\tag{\sf ASP-1} \\
    \mbox{subject to} &{\quad \widehat{\psi}_i(\br) \leq c_i + \pmb{\frac{\eta}{K}}  }, ~~~\forall i\in [m]\nonumber \\
    &{\quad \bB \br \leq \bd  }\nonumber \\
    &\quad l_j\leq r_j \leq u_j,~~\forall j\in [n]\nonumber
\end{align}
The following theorem establishes  upper bounds for both the feasibility and the objective value provided by a solution obtained from solving \eqref{prob:Approx-SP-2}. It generally shows that solving \eqref{prob:Approx-SP-2} will always yield a solution that is \(\cO(1/K)\) feasible to the original problem and yields an objective value within an \(\cO(1/K)\) neighborhood of the optimal value of [\textbf{SP}]. Consequently, solutions to \eqref{prob:Approx-SP-2} will converge linearly to an optimal solution to [\textbf{SP}] as the number of breakpoints \(K\) increases.
\begin{theorem}\label{th:overal-bound}
     Let $\widehat{\br}$ and $\br^*$ be optimal solutions to the approximate problem \eqref{prob:Approx-SP-2}  and the original one \textbf{[SP]}, we always have that $\widehat{\br}$ is $\frac{2\eta}{K}$-feasible to \textbf{[SP]}, and $|F(\br^*) - F(\widehat{\br})| \leq \frac{2\omega}{K}$, where $\eta$ and $\omega$ is defined in Lemmas \ref{lm:bound-F-hatF} and \ref{lm:bound-psi}.
\end{theorem}





\section{Constrained Dynamic Pricing}\label{sec:DP}

As discussed earlier, the commonly employed dynamic programming formulation for network RM suffers from the curse of dimensionality, as the state space expands exponentially with the number of resources. To address this issue, prior works typically rely on deterministic approximations or dynamic programming approximation approaches to break down the network problem into a set of single-resource problems \citep{liu2008choice,gallego1997multiproduct,zhang2017dynamic}. In the following, we briefly describe the deterministic approximation dynamic programming decomposition. More details can be found in the appendix.

The use of a \textbf{deterministic}  \textbf{approximation} model has been a popular approach for addressing the dynamic pricing problem. In our context, this can be achieved by replacing the probabilistic and discrete customer arrivals with a continuous fluid model characterized by an arrival rate \(\lambda\). Given a price vector \(\mathbf{r}\), the fraction of customers purchasing product \(j\) is denoted by \(P_j(\mathbf{r})\). Consequently, we let \(d = \lambda T\) represent the total customer arrivals over the time horizon \([0, T]\) and formulate the deterministic model as follows:
\begin{align}
 [\textbf{SP}^*] ~~~~~ \max_{\substack{\br\\l_j\leq r_j \leq u_j,~\forall j\in [n]}} &\left\{ F(\br) =  d\sum_{j=1}^n r_j P_j(\br)  \right\} \label{prob:SP-3} \\
    \mbox{subject to} &{\quad d\bA \bP(\br) \leq \bc  }\label{ctr:sp*-resource} \\
    &{\quad \bB \br \leq \bd  }\nonumber 
\end{align}
It can be observed that [\textbf{SP}$^*$] is a static pricing problem, for which the piece-wise linear approximation approach presented above can be directly utilized to achieve near-optimal solutions.

An alternative approach to deal with the curse of dimensionality is to employ a \textbf{dynamic decomposition approach }\citep{zhang2013assessing}. Specifically, for each fixed \( i \), the value function \( v_t(\bx) \) can be approximately decomposed by
$v_t(\bx) \approx v_{t,i}(x_i)  + \sum_{k \neq i} x_k \pi_k$,
for each \( t \) and \( x \), where \( \pi_k \) are the Lagrangian multipliers associated with constraints \eqref{ctr:sp*-resource}. The value \( v_t(\bx) \) is approximated by the sum of a nonlinear term of resource \( i \) and linear terms of all other resources. The Lagrangian multipliers \( \pi \) corresponding to constraint \eqref{ctr:sp*-resource} can be interpreted as the marginal value of an additional unit of each resource and incorporated into a dynamic programming decomposition approach \cite{zhang2013assessing}. With this decomposition, the dynamic program can be broken into a sequence of \( m \) one-dimensional dynamic programs where the value functions can be computed by solving static pricing problems of the form [\textbf{SP}]. With pricing constraints, our discretization techniques described above provide near-optimal solutions.

\section{Numerical Study}\label{sec:Numerical Results}
In this section, we provide numerical comparisons and analyses for both static and dynamic pricing under complex constraints. We will first describe some standard baselines, followed by numerical results for the static and dynamic pricing problems. 

\subsection{Baselines}

\subsubsection{Static Pricing:}
For solving [\textbf{SP}], we will compare our method described in Section \ref{sec:SP} with some standard baselines in the literature. In the following, we briefly discuss the baseline methods included in the experiments. More details can be found in the appendix.
\begin{itemize}
    \item \texttt{SP-DMIP:} This is our method, standing for Discretization and MILP. 
    \item \texttt{SP-SLSQP:}We use the \textit{Sequential Least Squares Programming} (SLSP) method \citep{kraft1988software} to solve the \textbf{[SP]}. 
    \item \texttt{SP-Trans:} Prior pricing methods typically transform the pricing problem into a convex optimization program with constraints on the purchasing probabilities. We adapted this approach for comparison with our method. Specifically, we first removed the pricing constraints (\(\bB\br\leq \bd\) and \(r_j \in [l_j, u_j], \forall j \in [n]\)) and converted [\textbf{SP}] into the purchasing-probability space, solving it using convex optimization. The resulting price solutions were then projected onto the feasible set to obtain feasible pricing solutions.
    \item \texttt{Gurobi:}  In this approach we use GUROBI's piece-wise linear approximation solvers \citep{GUROBI}  to approximately solve each step of the bisection described in Section \ref{sec:SP}.
\end{itemize}


\subsubsection{{Dynamic Pricing:}}
We compare our approach, which solves each step of the DPD using our discretization and MILP approximation method, with several baselines, similar to the static pricing case. These methods are briefly described below, with more details provided in the appendix.
\begin{itemize}
    \item \texttt{DP-DMIP}: This is our method, solving each step of the dynamic programming decomposition (DPD) using the discretization techniques described in Section \ref{sec:DP}.
    \item \texttt{DP-SLSQP}: This approach relies on the use of  \textit{SLSP} to directly solve each step of the DPD.
    \item \texttt{DP-Trans}: We use the same method as \texttt{SP-Trans} for solving each step of the DPD.
\end{itemize}

For dynamic pricing, we do not include GUROBI's PWLA in the comparison. The reason is that our experiments for \textbf{[SP]} below will show that GUROBI's PWLA is much slower than our approach (while providing the same performance guarantees). Therefore, for dynamic pricing problems, where static pricing problems need to be solved repeatedly, GUROBI's PWLA becomes inefficient. 

As mentioned earlier, one can solve the deterministic approximation [\textbf{SP}$^*$] and use the static solution as a fixed pricing strategy for the entire selling horizon. For the sake of comparison, we include this approach in our analysis. Specifically, we solve [\textbf{SP}$^*$] using three methods for static pricing: \texttt{SP-DMIP}, \texttt{SP-SLSQP}, and \texttt{SP-Trans}. The static solutions obtained are then compared with those derived from solving the dynamic programming decomposition in a simulation experiment.

All the experiments were implemented using Python and run on a server equipped with an AMD EPYC 7763 64-Core Processor. The MILPs were solved using GUROBI with default settings. We also utilized Scipy's nonlinear solvers for the approaches involving nonlinear optimization, specifically \texttt{SP-SLSQP}, \texttt{DP-SLSQP}, \texttt{SP-Trans}, and \texttt{DP-Trans}. We chose \( K=15 \) as the number of breakpoints used in our discretization approach. Experiments supporting this decision can be found in the appendix.

\subsection{Instance Generation}
Following prior works, for each set of instances, we randomly generated $T$ streams of demand arrivals, where the arrival in each period can be represented by a uniform$[0,1]$ random variable $X$. To generate product selection according to the choice probabilities \( P_j(\mathbf{r}) \), given the product prices \( \mathbf{r} \), an incoming customer selects product \( j \) if \( \sum_{k=1}^{j-1}P_k(\mathbf{r}) \leq X \leq \sum_{k=1}^{j}P_k(\mathbf{r}) \). We applied this procedure to create 20 different problem instances, which are labeled as cases 1 through 20.

The MNL choice parameters are generated as follows. The $a_j$ values for products are generated from a uniform $[10,100]$ distribution. The value $b_0$ is generated from a uniform $[0,20]$ and $b_j$ is generated from a uniform $[0,100]$.The lower bound $l_j$ and upper bound $u_j$ of price are generated by a uniform $[100,150]$ and $[250,400]$ distribution. The element of $\bA$ is generated from a uniform discrete distribution $[0,1]$. The parameters $\mathbf B$ and $ \mathbf d$  for the constraint $ \mathbf{B} \mathbf{r} \leq \mathbf{d} $ are generated as follows: For each instance set, we generate 3 additional constraints. The matrix $\bB$, with dimensions $3 \times n$ , is constructed such that each element is either 0 or 1, determined by a uniform discrete distribution $[0,1]$. Each row of $\mathbf{B}$ contains a specific number of value 1, given by a uniform discrete distribution $[0.5 n, 0.7n]$. This ensures variability in the constraints by introducing a controlled number of active entries per row. The vector $\mathbf{d} $, with dimensions $3 \times 1$, is calculated such that each element $d_i$ is given by a uniform distribution $[0.3 \sum_{j=1}^n u_j,0.5 \sum_{j=1}^n u_j] $. The arrival probability $\lambda$ is taken to be 1 in the selling horizon.

\subsection{Experiment Results}

\subsubsection{Comparison Results}

We conducted experiments on four scenarios:
(a) Small-size scenario with 2 resources, 3 products, and 200 time periods, each having a capacity of 60,   (b) medium-size scenario with 4 resources, 8 products, and 50 time periods, each with a capacity of 30,
  (c) large-size scenario with 16 resources, 80 products, and 50 time periods, each with a capacity of 30, and
  (d) the largest-size scenario with 16 resources, 80 products, and 200 time periods, each with a capacity of 120.
We first present a comparison of different methods for the constrained static pricing problem [\textbf{SP}]. Table \ref{table: revenue and running time of static pricing} reports the final revenues and running times for four methods: \texttt{SP-DMIP}, \texttt{SP-SLSQP}, \texttt{SP-Trans}, and GUROBI's PWLA. The results show that our method, \texttt{SP-DMIP}, consistently achieves better objective values than the others. While both \texttt{SP-DMIP} and \texttt{Gurobi} deliver competitive revenues, \texttt{SP-DMIP} has considerably shorter running times for large instances. This efficiency stems from \texttt{SP-DMIP}'s need for only one set of binary variables to approximate \(f_j(r_j)\) and \(g_j(r_j)\) for each product \(j\), unlike \texttt{Gurobi}, which requires a separate set of binary variables for each function. Consequently, \texttt{SP-DMIP} requires fewer additional binary variables, resulting in faster performance for large instances. Moreover, the ability to relax part of the binary variables in \texttt{SP-DMIP} (Theorem \ref{th:relax-MILP}) further enhances the method's efficiency. The other methods (i.e. \texttt{SP-Trans}  and \texttt{SP-SLSQP}) are fast but give significantly lower revenues, compared to \texttt{SP-DMIP} and \texttt{Gurobi}.

\begin{table}[htbp]
    \centering
    \caption{Revenue \& running time (seconds) comparison for static pricing}
    \label{table: revenue and running time of static pricing}
    \begin{tabular}{llrrrr}
        \toprule
        \textbf{Instance set}&  & \texttt{SP-DMIP} & \texttt{SP-SLSQP} & \texttt{SP-Trans} & \texttt{Gurobi} \\
        \midrule     
        \multirow{2}*{(2,3,200)} & Revenue &  \textbf{8364.85} &  5437.39 & 5437.39 & \textbf{8364.85} \\
        ~ & Time & 0.16 &  \textbf{0.01} & \textbf{0.01} & 0.24 \\

        \multirow{2}*{(4,8,50)} & Revenue & \textbf{8875.99} & 7373.31 & 7300.73& \textbf{8875.99} \\
        ~ & Time & 0.2 & \textbf{0.01} & 0.03 & 0.35 \\

        \multirow{2}*{(16,80,50)} & Revenue & \textbf{17723.31}  & 14451.38 &12302.23 & 17706.07 \\
        ~ & Time & 1.64 & 5.06 & \textbf{1.03} &  180.99 \\

        \multirow{2}*{(16,80,200)} & Revenue & \textbf{57779.33}  & 57430.93 & 30335.24 & 57567.08 \\
        ~ & Time & 1.97 & 7.97 & \textbf{0.52} & 182.17 \\
  
        \bottomrule
    \end{tabular}
\end{table}



\begin{table}
\centering
\caption{Objective values (i.e., revenues)  and running time (seconds) given by different methods for dynamic pricing}
\label{table: bounds}
\begin{tabular}{c|l|ccc|ccc} 
\toprule
\multicolumn{1}{c}{}                      &         & \multicolumn{3}{c|}{\texttt{DP-}}                & \multicolumn{3}{c}{\texttt{SP-}}                  \\ 
\midrule
\multicolumn{1}{c}{\textbf{Instance set}} &         & \texttt{DMIP}& \texttt{SLSQP}& \texttt{Trans}& \texttt{DMIP}& \texttt{SLSQP}& \texttt{Trans}\\ 
\midrule
\multirow{2}{*}{(2,3,200)}                & Revenue & \textbf{3802.99 } & 2525.17  & 2559.65  & \textbf{3839.44}   & 2525.17  & 1117.87  \\
                                          & Time    & 2777.45           & 17.41    & 146.61   & 0.49               & 0.0031   & 0.07     \\ 
\hline
\multirow{2}{*}{(4,8,50)}                 & Revenue & \textbf{2616.79}  & 2494.02  & 1337.62  & \textbf{2645.1}    & 2507.51  & 1334.88  \\
                                          & Time    & 365.99            & 13.95    & 16.75    & 0.98               & 0.021    & 0.04     \\ 
\hline
\multirow{2}{*}{(16,80,50)}               & Revenue & \textbf{8785.89}  & 8520.8   & 7233.52  & \textbf{8791.37 }  & 8593.05  & 5304.42  \\
                                          & Time    & 2608.53           & 808.31   & 2693.78  & 3.46               & 0.06     & 0.46     \\ 
\hline
\multirow{2}{*}{(16,80,200)}              & Revenue & \textbf{34133.28} & 26945.74 & 30499.29 & \textbf{34133.36 } & 33764.44 & 32492.6  \\
                                          & Time    & 18613.47          & 11246.01 & 26899.98 & 2.16               & 0.1      & 2.3      \\
\bottomrule
\end{tabular}
\end{table}

In the subsequent part of the section, we provide  a comparison for the dynamic pricing problem. Due to the prohibitively slow performance of \texttt{Gurobi} for large instances in the static setting (shown in Table \ref{table: revenue and running time of static pricing}), we have excluded it from this comparison. Table \ref{table: bounds} presents the objective values (i.e., revenues) achieved by various methods when solving the dynamic pricing instances using the  dynamic programming  or solving the deterministic approximation decomposition [\textbf{SP}$^*$] (as described in Section \ref{sec:DP}). Here we note that, for static pricing with constraints on both prices and purchasing probability, \texttt{SP-DMIP} can guarantee near-optimal solutions, whereas \texttt{SP-Trans} and \texttt{SP-SLSQP} yield suboptimal solutions. The results in Table \ref{table: bounds} clearly demonstrate that \texttt{SP-DMIP} generates significantly higher revenues compared to other static pricing approaches. Furthermore, for the dynamic programming decomposition approaches, \texttt{DP-DMIP} also achieves higher revenues compared to other dynamic  methods (\texttt{DP-SLSQP} and \texttt{DP-Trans}). These results confirm the advantages of our piecewise linear approximation approach for effectively addressing both static and dynamic pricing problems with pricing constraints. It is also important to note that the revenues generated by dynamic pricing approaches consistently fall below those achieved by the deterministic approximation. This observation  well aligns the findings in previous studies \citep{zhang2013assessing} which state that the approximate dynamic decomposition yields tighter upper bounds compared to the static approximation.

In terms of computing times, all methods perform quickly for the static methods. For dynamic programming decomposition, \texttt{DP-DMIP} requires more time than \texttt{DP-SLSQP}, which is understandable since \texttt{DP-DMIP} involves solving several MILPs, whereas \texttt{DP-SLSQP} only needs to solve nonlinear (continuous) problems, typically a faster process. Moreover, \texttt{DP-DMIP} is slower than \texttt{DP-Trans} for small-sized instances but is faster for larger-sized instances. Overall, our numerical comparisons indicate that \texttt{DP-DMIP} delivers significantly better objective values than other methods while maintaining reasonable computing times.

\begin{table}
\centering
\caption{Simulation results (i.e., revenues) given by different methods.}
\label{table: simulation results}
\begin{tabular}{c|l|ccc|ccc} 
\toprule
\multicolumn{1}{c}{}                      &      & \multicolumn{3}{c|}{\texttt{DP-}}                  & \multicolumn{3}{c}{\tt{SP-}}          \\ 
\hline
\multicolumn{1}{c}{\textbf{Instance set}} &      & \tt{DMIP}               & \tt{SLSQP}    & \tt{Trans}     & \tt{DMIP}      & \tt{SLSQP}    & \tt{Trans}     \\ 
\hline
\multirow{2}{*}{(2,3,200)}                & mean & \textbf{3771.1}    & 2505.0   & 2578.68   & 3658.18   & 2505.0   & 11118.90  \\
                                          & std  & \textbf{(159.25)}  & (409.72) & (389.03)  & (218.92)  & (409.72) & (362.22)  \\ 
\hline
\multirow{2}{*}{(4,8,50)}                 & mean & \textbf{2613.0}    & 2542.0   & 1280.06   & 2574.49   & 2522.0   & 1278.50   \\
                                          & std  & \textbf{(175.23)}  & (197.45) & (347.32)  & (154.28)  & (177.22) & (345.29)  \\ 
\hline
\multirow{2}{*}{(16,80,50)}               & mean & \textbf{8987.97}   & 7748.64  & 8070.49   & 8931.4    & 8584.51  & 5079.66   \\
                                          & std  & \textbf{(734.84)}  & (670.49) & (676.33)  & (620.35)  & (493.51) & (162.26)  \\ 
\hline
\multirow{2}{*}{(16,80,200)}              & mean & \textbf{34043.92}  & 23006.52 & 29957.40  & 34042.59  & 33657.30 & 32110.40  \\
                                          & std  & \textbf{(1405.86)} & (870.21) & (1735.82) & (1402.12) & (784.52) & (970.61)  \\
\bottomrule
\end{tabular}
\end{table}

Table \ref{table: bounds} presents only the objective values of the approximate problems (i.e., dynamic programming decomposition and static approximation), which do not reflect how the resulting pricing solutions perform in the actual dynamic pricing context. To further evaluate the performance of the solutions derived for the original dynamic pricing problem, we conducted a simulation experiment. Specifically, we repeatedly simulated customer choice behavior and collected revenue data across the selling horizons for each pricing solution. This process was repeated 20 times for each pricing solution, and the average revenue obtained offers an approximation of the actual expected revenue, which will converge to the true value as the number of repetitions increases. Table \ref{table: simulation results} reports the means and standard deviations of the revenues from different pricing solutions.
The simulation results clearly demonstrate that our approach, \texttt{DP-DMIP}, consistently delivers higher revenue compared to all other methods. Interestingly, \texttt{SP-DMIP}, despite only solving the static approximation of the dynamic pricing problem, provides competitive revenues compared to \texttt{DP-DMIP} and outperforms other methods. Moreover, both \tt{DP-DMIP} and \tt{SP-DMIP} yield relatively smaller standard deviations compared to the other approaches.

In summary, our numerical comparisons demonstrate that, with constraints on both prices and purchasing probabilities, our piece-wise linear approximation approaches consistently outperform other approaches based on nonlinear optimization solvers or GUROBI's PWLA. In particular, while \texttt{DP-DMIP} provides better revenues, its deterministic counterpart \texttt{SP-DMIP} also offers competitive revenues in both estimation and simulation tests. These results highlight the advantages of our approaches in solving constrained pricing problems in both static and dynamic settings.

\section{Conclusion}\label{sec:concl}

In this study, we have presented a general approach to tackling pricing problems with constraints, which are known to be highly nonconvex and challenging to solve to optimality. We explored a piece-wise linear approximation method that allows us to approximately solve the static pricing problem to any desired precision by solving a sequence of MILPs. We further established bounds on the approximation errors resulting from solving the approximate problem and demonstrated how discretization techniques can be used for constrained dynamic pricing. Extensive experiments show that our methods outperform other standard baselines, consistently yielding higher revenue in our simulation tests. Future work should focus on exploring more advanced customer choice models, such as nested or cross-nested logit models \citep{train2009discrete}, or delving into the structure of dynamic pricing problems to develop better decomposition methods.

\bibliographystyle{informs2014} 
\bibliography{Reference.bib} 

\begin{APPENDICES}
\clearpage
\begin{center}
    {\huge Appendix}
\end{center}
\section{Missing Proofs}\label{appdx:proofs}

\subsection{Proof of Proposition \ref{prop:bisection-gap}}
\begin{proofs}
    The  proposition can be verified quite straightforwardly,  as we can see that  after N steps, the length of the  searching interval becomes $(U-L)/2^N$. This directly implies that the algorithm will stop after $N \leq \ceil{\log\left(\frac{U-L}{\xi}\right)}$ steps, noting that  when $N =  \leq {\log\left(\frac{U-L}{\xi}\right)}$ we have $(U-L)/2^N = \xi$.
    Now, let  $L^*,U^*$ be the lower bound and upper bound after N steps of the binary search,  we have $U^*-L^*  = (U-L)/2^N$, and 
    \begin{align}
        \max\limits_{\br\in\widehat{\cR}}~ \left\{ \sum_{j\in [n]} \widehat{f}_j(r_j) - L^* \left(1+\sum_{j\in[n]} \widehat{g}_j(r_j)\right) \right\} \geq 0 \nonumber\\
              \max\limits_{\br\in\widehat{\cR}}~ \left\{ \sum_{j\in [n]} \widehat{f}_j(r_j) - U^* \left(1+\sum_{j\in[n]} \widehat{g}_j(r_j)\right) \right\} \leq 0 \nonumber
    \end{align}
So, if  $\widehat{\br}$ is  an optimal solution to    $ \max\limits_{\br\in\widehat{\cR}}~ \left\{ \sum_{j\in [n]} \widehat{f}_j(r_j) - L^* \left(1+\sum_{j\in[n]} \widehat{g}_j(r_j)\right) \right\} $, we have 
\begin{align*}
    \widehat{F}(\widehat{\br}) \geq L^*;~~~\text{ }
    \max_{\br\in \widehat{\cR}} \widehat{F}(\br) \leq  U^* 
\end{align*}
Which implies
\[
\widehat{F}(\br^*) - \widehat{F}(\widehat{\br}) \leq U^*-L^* = \xi,
\]
as desired. 
\end{proofs}

\subsection{Proof of Lemma \ref{lm:1}}
\begin{proofs}
    It is obvious to see that $g_j(r_j) = e^{(a_j-r_j)/b_j}$ is convex in $r_j$. For $f_j(r_j) = r_je^{(a_j-r_j)/b_j}$, we take it's first and second-order derivatives  and write
\begin{align*}
    \frac{\partial f_j(r_j)}{\partial r_j} &= e^{(a_j-r_j)/b_j} - \frac{r_j}{b_j} e^{{(a_j-r_j)/b_j}} \\
    \frac{\partial^2 f_j(r_j)}{\partial r_j\partial r_j} &= -\frac{2}{b_j}e^{(a_j-r_j)/b_j} + \frac{r_j}{b^2_j} e^{{(a_j-r_j)/b_j}}
\end{align*}
So it can be seem that $ \frac{\partial^2 f_j(r_j)}{\partial r_j\partial r_j} \geq 0$ when $r_j \geq 2b_j$ and $ \frac{\partial^2 f_j(r_j)}{\partial r_j\partial r_j} \leq  0$ otherwise, validating the claim
\end{proofs}

\subsection{Proof of Lemma \ref{lm:2}}

\begin{proofs}
   For any $k\leq K_j-2$, recall that 
   \[
    \gamma^g_{jk} = \frac{g_j((k+1)\Delta_j) - g_j(k\Delta_j)}{\Delta_j};  ~ \gamma^g_{j,k+1} = \frac{g_j((k+2)\Delta_j) - g_j((k+1)\Delta_j)}{\Delta_j}.
   \]
The Mean Value Theorem implies that there is $c^g_k\in [k\Delta_j,(k+1)\Delta_j]$  and $c^g_{k+1}\in [(k+1)\Delta_j,(k+2)\Delta_j]$ such that $\gamma^g_{jk} = g'(c^g_k)$  and $\gamma^g_{j,k+1} = g'(c^g_{k+1})$. Moreover, since $g_j(r_j)$ is strictly convex in $r_j$, $g'_j(r_j)$ is strictly  increasing in $r_j\in [l_j,u_j]$. Thus,  $ g'(c^g_k)<  g'(c^g_{k+1})$, implying $\gamma^g_{jk}< \gamma^g_{j,k+1}$, as desired.

Similarly, from Lemma \ref{lm:1} we know that $f_j(r_j)$ is convex when  $r_j\leq 2b_j$. Thus, for any $ k\leq \lfloor\frac{2b_j}{\Delta_j}\rfloor -2 $ we have $(k+2)\Delta_j \leq 2b_j$. Moreover, from the Mean Value Theorem, for any $k$ such that   $(k+2)\Delta_j \leq b_j$, there are  $c^f_k\in [k\Delta_j,(k+1)\Delta_j]$  and $c^f_{k+1}\in [(k+1)\Delta_j,(k+2)\Delta_j]$ such that $\gamma^f_{jk} = f'(c^f_k)$  and $\gamma^f_{j,k+1} = f'(c^f_{k+1})$. Since  $f(r_j)$ is convex for $r_j\leq 2b_j$, $f'(r_j)$ is monotonically increasing in $r_j$, implying $f'(c^f_k)\leq f'(c^f_{k+1})$, or $\gamma^f_{jk} \leq \gamma^f_{j,k+1}$ as desired.

For the case that $k\geq\lceil\frac{2b_j}{\Delta_j}\rceil $,  in a similar way, we see  that the function is concave when $r_j \geq 2b_j$. We can also use the Mean Value Theorem to see that $\gamma^f_{jk} \geq \gamma^f_{j,k+1}$. We complete the proof.
\end{proofs}

\subsection{Proof of Lemma \ref{lm:bound-F-hatF}}
\begin{proofs}
 We first see that, for any $r_j,r'_j \in [l_j,u_j]$ we have 
  \begin{align*}
   |f_j(r_j) - f_j(r'_j)| &\le \max_{c\in [l_j,u_j]} f'_j(c) |r_j-r_j'|\\
        & = \max_{c\in [l_j,u_j]} \left|e^{(a_j-c)/b_j} - \frac{c}{b_j} e^{(a_j-c)/b_j} \right| |r_j-r_j'|\\
        &\leq  e^{(a_j-l_j)/b_j}\left(1+\frac{u_j}{b_j}\right) |r_j-r'_j|
  \end{align*}
  and
\begin{align}
     |g_j(r_j) - g_j(r'_j)|&\leq \max_c |g_j'(c)| |r_j-r'_j|\nonumber\\
     & \max_c |\frac{1}{b_j} e^{(a_j-c)/b_j}| |r_j-r'_j|\nonumber\\
     &\le \frac{e^{(a_j-l_j)/b_j}}{b_j} |r_j-r'_j|\label{eq:Lipschitz-g}
\end{align}
We now bound the gaps  $|f_j(r_j)-\widehat{f}_j(r_j)|$  and $|g_j(r_j)-\widehat{g}_j(r_j)|$.  Let $k = \lfloor \frac{r_j}{\Delta_j}\rfloor$ and $\psi_j = r_j - k\Delta_j$.  Without  the loss of generality, assume that $f_j(k\Delta_j)\ge f_j(r_j)((k+1)\Delta_j)$ (the other case can be handled in a similar way). We then compute $\widehat{f}_j(r_j)$ as 
    \[
     \widehat{f}_j(r_j) =  f_j(k\Delta_j) + \frac{f_j((k+1)\Delta_j- f_j(k\Delta_j))}{\Delta_j} \psi_j
  \]
    We then see that $\widehat{f}_j(r_j) \leq f_j(k\Delta_j)$ and 
    \[
      \widehat{f}_j(r_j) \geq  f_j(k\Delta_j) + \frac{f_j((k+1)\Delta_j- f_j(k\Delta_j))}{\Delta_j} \Delta_j = f_j((k+1)\Delta_j)
    \]
Thus $f_j(k\Delta_j)\geq \widehat{f}_j(r_j)\geq f_j((k+1)\Delta_j)$, implying
    \begin{align}
  |f_j(r_j) - \widehat{f}_j(r_j)| &\le \max\left\{|f_j(r_j) - f( k \Delta_j)|; ~ |f_j(r_j) - f( (k+1) \Delta_j)|\right\}  \\
  &\le e^{(a_j-l_j)/b_j}\left(1+\frac{u_j}{b_j}\right)  \Delta_j \label{eq:bound-f}
  \end{align}
Similarly for the gap  $|g_j(r_j)-\widehat{g}_j(r_j)|$  we also have
\begin{equation}\label{eq:bound-g}
 |g_j(r_j) - \widehat{g}_j(r_j)| \leq \frac{e^{(a_j-l_j)/b_j}}{b_j}  \Delta_j  
\end{equation}

We further remark that $g_j(r_j)\leq \widehat{g}_j(r_j)$ for all $r_j\in [l_j,u_j]$, which is due to the convexity of $g_j(\cdot)$. 
We are now ready to evaluate the bound $|F(\br) - \widehat{F}(\br)|$. 
For ease of notation, let $U(\br) = \sum_{j\in [n]} f_j(r_j)$, $V(\br) = \sum_{j\in [n]} g_j(r_j)$, $\widehat{U}(\br) = \sum_{j\in [n]} \widehat{f}_j(r_j)$, $\widehat{V}(\br) = \sum_{j\in [n]} \widehat{g}_j(r_j)$. 
From the bounds in \eqref{eq:bound-f} and \eqref{eq:bound-g} we have 
\begin{align*}
    |U(\br) - \widehat{U}(\br)| &\leq \frac{1}{K}\sum_{j\in [n]}  e^{(a_j-l_j)/b_j}\left(1+\frac{u_j}{b_j}\right)(u_i-l_i)\\
       |V(\br) - \widehat{V}(\br)| &\leq \frac{1}{K}\sum_{j\in [n]}  \frac{e^{(a_j-l_j)/b_j}}{b_j}(u_i-l_i)
\end{align*}
We further see that $\widehat{V}(\br)\ge V(\br)$. For the sake of simplicity, let us denote $\alpha = \sum_{j\in [n]}  e^{(a_j-l_j)/b_j}\left(1+\frac{u_j}{b_j}\right)(u_i-l_i)$ and $\beta= \sum_{j\in [n]}  \frac{e^{(a_j-l_j)/b_j}}{b_j}(u_i-l_i)$. To bound the gap $|F(\br) - \widehat{F}(\br)|$, we consider the following cases
\begin{itemize}
    \item If $F(\br)\geq \widehat{F}(\br)$: We first remark that 
         \[
         F(\br) = \frac{\sum_{j\in  [n]} r_j g_j(r_j)}{1+\sum_{j\in [n]}g_j(r_j)} \leq \max_{j\in [n]} \{r_j\} \stackrel{def}{=} r^u
        \]
        Therefore,
        \[
           \frac{U(\br)}{V(\br)} \leq \frac{U(\br) + (\widehat{V}(\br) - V(\br)) r^u}{V(\br) + (\widehat{V}(\br) - V(\br))}  = \frac{U(\br) + (\widehat{V}(\br) - V(\br)) r^u}{\widehat{V}(\br) }
        \]
 So we have
\begin{align*}
    F(\br) - \widehat{F}(\br) &\leq \frac{U(\br) - \widehat{U}(\br) + (\widehat{V}(\br) - V(\br)) r^u}{\widehat{V}(\br)} \\
    &\leq \left(\frac{\alpha}{K}+ \frac{\beta}{K}r^u\right) \frac{1}{\min_{\br} V(\br)}\\
    & \leq \left(\frac{\alpha}{K}+ \frac{\beta}{K}r^u\right) \frac{1}{1+\sum_{j\in [n]} e^{(a_j-u_j)/b_j}}\\
\end{align*}
\item If $F(\br)< \widehat{F}(\br)$, we write 
f\begin{align*}
    \widehat{F}(\br) - F(\br)&=\frac{\widehat{U}(\br)}{\widehat{V}(\br)} - \frac{U(\br)}{V(\br)}\\
     &\leq \frac{\widehat{U}(\br) - U(\br)}{{V}(\br)} \\
     &\leq \frac{\alpha}{K} \frac{1}{1+\sum_{j\in [n]} e^{(a_j-u_j)/b_j}}
\end{align*}
\end{itemize}
Combining the two cases, we get
\[
  | \widehat{F}(\br) - F(\br)|\leq \frac{1}{K}\frac{\alpha +\beta r^u}{1+\sum_{j\in [n]} e^{(a_j-u_j)/b_j}} = \frac{\omega}{K}
\]
which completes the proof.
\end{proofs}

\subsection{Proof of Theorem \ref{th:relax-MILP}}

\begin{proofs}
To prove the equivalence, let \((\bw^*, \bz^*, \br^*)\) be an optimal solution to \eqref{prob:MILP-2} (i.e., the MILP approximation where part of the binary variables are relaxed). We need to prove that this solution is also equivalent to an optimal solution to \eqref{prob:MILP}. We will consider two cases: \(K_j > \lfloor\frac{2b_j}{\Delta_j}\rfloor\) or \(K_j \leq \lfloor\frac{2b_j}{\Delta_j}\rfloor\).

When \(K_j > \lfloor\frac{2b_j}{\Delta_j}\rfloor\), we leverage Lemma \ref{lm:2} to see that for all \(k \leq \lfloor\frac{2b_j}{\Delta_j}\rfloor - 1\) we have \(\gamma^g_{jk} \geq \gamma^g_{j,k-1}\) and \(\gamma^f_{jk} > \gamma^f_{j,k-1}\). So, if we let \(\tau_{jk} = \gamma^f_{jk} - \delta \gamma^g_{jk}\) (for ease of notation), then we see that \(\tau_{jk}\) are strictly increasing in \(k\), for all \(k \leq \lfloor\frac{2b_j}{\Delta_j}\rfloor - 1\), i.e.,
\[ \tau_{j1} < \tau_{j2} < \ldots < \tau_{j,\lfloor\frac{2b_j}{\Delta_j}\rfloor - 1} .\]
Moreover, the objective function in \eqref{prob:MILP-2} involves the term 
\[ \sum_{j \in [n]} \sum_{k = 0}^{\lfloor\frac{2b_j}{\Delta_j}\rfloor - 1} \Delta_j \tau_{jk} w_{jk}. \]
where \(w_{jk}\), \(j \in [n]\), \(k \in \{0, \ldots, K_j - 1\}\), satisfy the following constraints:
\begin{align}
    & z_{jk} \geq z_{j,k+1} \quad k = 0, \ldots, K-2, j \in [n] \label{ctr:-1} \\
    & z_{jk} \leq w_{jk} \leq 1 \quad k = 0, \ldots, K-1, j \in [n] \label{ctr:-2} \\
    & w_{j,k+1} \leq z_{j,k} \quad k = 0, \ldots, K-2, j \in [n] \label{ctr:-3} \\
    & r_j = \Delta_j \sum_{k = 0}^{K-1} w_{jk} \quad j \in [n] \label{ctr:-4} \\
    & \bB \br \leq \bd \label{ctr:-5} \\
    & z_{jk} \in [0,1], 0 \leq k \leq \left\lfloor \frac{2b_j}{\Delta_j} \right\rfloor - 1, \text{ and } z_{jk} \in \{0,1\}, k \geq \left\lfloor \frac{2b_j}{\Delta_j} \right\rfloor .\label{ctr:-6}
\end{align}

So, if there is any \(k \geq \left\lfloor \frac{2b_j}{\Delta_j} \right\rfloor\) such that \(z^*_{jk} = 1\), then we see that \(z^*_{jk} = 1\) for all \(k \leq \left\lfloor \frac{2b_j}{\Delta_j} \right\rfloor - 1\), implying that \((\bw^*, \bz^*, \br^*)\) is feasible to \eqref{prob:MILP}, thus is optimal to \eqref{prob:MILP}. Otherwise, if \(z^*_{jk} = 0\) for all \(k \geq \left\lfloor \frac{2b_j}{\Delta_j} \right\rfloor\), we see that, for any \(j \in [n]\), the monotonicity of \(\tau_{jk}\) implies that if there are \(k, h \in \{0, \ldots, \left\lfloor \frac{2b_j}{\Delta_j} \right\rfloor - 1\}\) such that \(k < h\) and \(w^*_{jk} < 1\) and \(w^*_{jh} > 0\), then we can always choose a better solution by increasing \(w^*_{jk}\) and reducing the value of \(w^*_{jh}\). Mathematically, there is \(\epsilon > 0\) such that we can always choose a solution \((\overline{\bw}, \overline{\bz}, \overline{\br})\) such that \(\overline{\br} = \br^*\), \(\overline{\bw} = \bw^*\) except \(\overline{w}_{jk} = w^*_{jk} + \epsilon\) and \(\overline{w}_{jh} = w^*_{jh} - \epsilon\), and \(\overline{\bz}\) is chosen to satisfy the constraints \eqref{ctr:-1} - \eqref{ctr:-6} (together with \(\overline{\bw}, \overline{\br}\)). We then see that \((\overline{\bw}, \overline{\bz}, \overline{\br})\) always yields a higher objective value than that given by \((\bw^*, \bz^*, \br^*)\). Thus, for any \(k, h \in \{0, \ldots, \left\lfloor \frac{2b_j}{\Delta_j} \right\rfloor - 1\}\) such that \(k < h\), we should either have \(w^*_{jk} = 1\) or \(w^*_{jh} = 0\). This directly implies that there is an index \(k^* \in \{0, \ldots, \left\lfloor \frac{2b_j}{\Delta_j} \right\rfloor - 1\}\) such that \(w^*_{jk} = 1\) for all \(k < k^*\), and \(w^*_{jk} = 0\) for all \(k > k^*\). We then can construct a binary solution \(\widehat{\bz}\) such that \(\widehat{z}_{jk} = 1\) when \(w^*_{jk} > 0\), and \(\widehat{z}_{jk} = 0\) otherwise and see that \((\overline{\bw}, \widehat{\bz}, \overline{\br})\) is optimal for the relaxed problem \eqref{prob:MILP-2}. Since \((\overline{\bw}, \widehat{\bz}, \overline{\br})\) is feasible to \eqref{prob:MILP}, this is also optimal to \eqref{prob:MILP}, confirming the equivalence.

For the second case when \(K_j \leq \lfloor\frac{2b_j}{\Delta_j}\rfloor\), we see that \(\tau_{jk}\) is strictly increasing for all \(k \in \{0, \ldots, K_j - 1\}\). We can use the same arguments as in the above case to see that one can construct a solution \((\overline{\bw}, \widehat{\bz}, \overline{\br})\) from \((\bw^*, \bz^*, \br^*)\) that is optimal to \eqref{prob:MILP}. This case happens  when $u_j\leq 2b_j$   and we can relax all the binary constraints in \eqref{prob:MILP}. This confirms the second statement of the theorem (i.e., all the binary variables in \eqref{prob:MILP} can be safely relaxed). This completes the proof.
\end{proofs}

\subsection{Proof of Theorem \ref{th:overal-bound}}
\begin{proofs}  
We first note that, since $\br^*$ is  feasible to \textbf{[SP]}, we should have  $\psi_i(\br) \leq c_i$, $\forall i\in [m]$. Thus
\[
\widehat{\psi}_i(\br^*) \stackrel{(a)}{\leq} \psi_i(\br^*) + \frac{\eta}{K} \leq c_i+\frac{\eta}{K}
\]
implying that $\br^*$ is feasible to \eqref{prob:Approx-SP-2}. It follows that $\widehat{F}(\br^*)\leq \widehat{F}(\widehat{\br})$.
Moreover, using the bound established in Lemma \eqref{lm:bound-psi}  we have 
\begin{align}
{\psi}_i(\widehat{\br}) &\leq \widehat{\psi}_i(\widehat{\br}) + \frac{\eta}{K}  \nonumber\\
&\leq c_i + \frac{2\eta}{K}, ~~\forall i\in [m] \nonumber
\end{align}
which implies that $\widehat{\br}$ is $\frac{2\eta}{K}$ feasible to \textbf{[SP]}. To bound the gap between $F(\br^*)$ and $F(\widehat{\br})$, we write
 \begin{align}
     F(\br^*) &\stackrel{(c)}{\leq} \widehat{F}(\br^*) + \frac{\omega}{K} 
     \stackrel{(d)}{\leq} \widehat{F}(\widehat{\br}) + \frac{\omega}{K}
     \stackrel{(e)}{\leq} {F}(\widehat{\br}) + \frac{2\omega}{K}\nonumber
 \end{align}
where $(c)$ is due to Lemma \ref{lm:bound-F-hatF}, $(d)$ is because  $\br^*$ is feasible to \eqref{prob:Approx-SP-2} and the last inequality $(e)$ is, again, due to Lemma \ref{lm:bound-F-hatF}. This directly implies that $|F(\br^*) - F(\widehat{\br})|\leq \frac{2\omega}{K}$ as desired. 

\end{proofs}

\section{Constrained Dynamic Pricing}\label{sec-appendix:DP}
\label{section DP}
In this section, we provide more details about the deterministic approximation and dynamic programming decomposition presented in Section \ref{sec:DP}.

The use of a deterministic and continuous approximation model has been a popular approach for addressing the dynamic pricing problem. 
It can be observed that [\textbf{SP}$^*$] is a static pricing problem, for which the piece-wise linear approximation approach presented above can be directly utilized to achieve near-optimal solutions. As discussed in prior work, a solution to [\textbf{SP}$^*$] can serve as a vector of prices. Since \(\mathbf{r}^*\) is a constant vector, independent of time or inventory levels, it results in a static pricing policy where prices remain fixed throughout the selling period. Additionally, without the price constraints \(\mathbf{B}\mathbf{r} \leq \mathbf{d}\), the Lagrangian multipliers \(\pi\) associated with constraints \eqref{ctr:sp*-resource} can be interpreted as  the value of an additional unit of each resource and can be employed in the dynamic programming decomposition discussed in the next section \citep{zhang2013assessing,zhang2009approximate}.

For the  dynamic programming decomposition approach to solve  [\textbf{DP}],  the value function $v_t(\bx)$ can be approximated by 
\begin{align}
\label{decomposition of value function}
v_t(\bx) \approx v_{t,i}(x_i)  + \sum_{k \neq i} x_k \pi_k 
\end{align}
for each i, $t$ and $x$. The value $v_t(\bx)$ is approximated by the sum of a nonlinear term of resource $i$ and linear terms of all other resources. The Lagrangian multipliers $\pi$ corresponding to constraint \eqref{ctr:sp*-resource} can be interpreted as the marginal value of an additional unit of each resource and incorporated into a dynamic programming decomposition approach \cite{zhang2013assessing}. Using (\ref{decomposition of value function}) in \textbf{[DP]} and adapt the feasible set $R_t(\bx)$ to $R_t(x,\bc_{-i})$, we obtain

\begin{align}
[\textbf{DP}_i]& \quad v_{t,i}(x_i) \label{value function i} \\
&=\max \limits_{r_t \in \cR_t(x_i,\bc_{-i})} \left\{ \sum_{j=1}^{n} \lambda P_j(r_t)[r_{t,j} -\sum_{k\neq i} a_{kj}\pi_k+ v_{t+1,i}(x_i-a_{ij}) - v_{t+1,i}(x_i)] \right\} + v_{t+1,i}(x_i) \nonumber\\
&= \max \limits_{r_t \in \cR_t(x_i,\bc_{-i})} \left\{\sum_{j=1}^{n}\lambda P_j(r_t)[r_{t,j} -\sum_{k\neq i} a_{kj}\pi_k-\triangle_{j}v_{t+1,i}(x_i)]\right\} + v_{t+1,i}(x_i) \nonumber
\end{align}

The boundary conditions are defined as $v_{T+1,i}(\bx) = 0$ for all values of $\bx$, and $v_{t,i}(\textbf{0}) = 0$ for all time steps $t$. In the above statement, $(x_i,\bc_{-i})$ represents an $m$-vector, where the $i$-th component is denoted as $x_i$ and the $k$-th component (for $k\neq i$) is denoted as $c_k$. By solving the set of $m$ one-dimensional dynamic programs, we can determine the values of $v_{t,i}(x_i)$ for each $i$.

Using similar techniques, the maximization in \textbf{[DP$_i$]} for each state $x_i$ and time $t$ can be reformulated as an nonlinear optimization problem similar to \textbf{[SP]}. At each time $t$, for each fixed resource $i$, solving (\ref{value function i}) is same as maximizing the first part of equation (\ref{value function i}). Therefore, we aim to solve the following problem at each stage $t$:

\begin{align}
    \xi_t(x_i) =  & \max \limits_{\br_t \in \cR_t(x_i,\bc_{-i})} \sum_{j=1}^n \lambda P_j(\br_t)\left[r_{t,j} -\sum_{k\neq i} a_{kj}\pi_k-\triangle_{j}v_{t+1,i}(x_i)\right] \label{dp_obj_1}\\
    \mbox{subject to} & \quad \quad l_j \leq r_{t,j} \leq u_j   \quad  \forall j \in [n] \label{dp_constraint_1} \\
& \quad \quad \lambda \bA \bP(\br_t) \leq (x_i, \bc_{-i}) \nonumber \\
& \quad \quad \bB \br_t \leq \bd\nonumber
\end{align}
This is similar to the static pricing problem [\textbf{SP}], with a slight change in the objective function. Therefore, the discretization approach used for solving [\textbf{SP}] can be employed to compute \(\xi_t(x_i)\). Specifically, we can use the Dinkelbach transform to convert \eqref{dp_obj_1} into a binary search procedure over a scalar \(\delta \in \mathbb{R}\), where, at each step, we need to solve the following MILP:

\begin{align}
[\textbf{MILP}^t_i] ~~~~~ \max\limits_{\bw,\bz,\br_t} & \sum_{j\in [n]} (f^t_j(l_j)-\delta g^t_j(l_j)) + \sum_{j \in [n]} \sum_{k = 0}^{K-1}\Delta^t_j(\gamma^{ft}_{jk}-\delta \gamma^{gt}_{jk})w_{jk}-\delta  \nonumber\\
\mbox{subject to} \quad &   z_{jk}\geq z_{j,k+1}  \quad k= 0,\ldots,K-2, j \in [n] \nonumber\\
& {z_{jk}} \leq w_{jk} \leq 1, \quad k = 0,\ldots,K-1, j \in [n]  \nonumber \\
& w_{j,k+1} \leq z_{j,k}  \quad k = 0,\ldots,K-2, j \in [n]  \nonumber\\
&r_{t,j} = \Delta_j\sum_{k=0}^{K-1} w_{jk},\quad \forall j\in [n]\nonumber\\
& \sum_{j\in [n]} (a_{ij} - x_i) \left(g^t_j(l_j) + \sum_{k=0}^{K-1} \gamma^{gt}_{jk} w_{jk}\right) \leq x_i \nonumber \\
& \sum_{j\in [n]} (a_{i'j} - c_{i'}) \left(g^t_j(l_j) + \sum_{k=0}^{K-1} \gamma^{gt}_{jk} w_{jk}\right) \leq c_{i'},~~\forall i'\in [m]\backslash \{i\}\nonumber \\
&\bB\br_t \leq \bd \nonumber\\
& z_{jk} \in \{0,1\}, w_{jk} \in \mathbb{R} \quad \forall k\in[K], j \in [n] \label{ctr-milp-8}
\end{align}

where
\begin{align*}
    f^t_j(r_{t,j}) &= \left(r_{t,j} -\sum_{k\neq i} a_{kj}\pi_k-\triangle_{j}v_{t+1,i}(x_i)\right)\exp\left(\frac{a^t_j-r_{t,j}}{b^t_j}\right),~\forall t\in [T], j\in [n] \\
    g^t_j(r_{t,j})&= \exp\left(\frac{a^t_j-r_{t,j}}{b^t_j}\right),~\forall t\in [T], j\in [n] \\
    \Delta^t_j &= \frac{u_j-l_j}{K}, ~\forall t\in [T], j\in [n]\\
    \gamma^{ft}_{jk} &= \frac{f_j^t((k+1)\Delta^t_j) - f_j^t(k\Delta^t_j)}{\Delta^t_j}, ~\forall t\in [T], j\in [n], k\in\{0,\ldots,K\}\\
    \gamma^{gt}_{jk} &= \frac{g_j^t((k+1)\Delta^t_j) - g^t_j(k\Delta^t_j)}{\Delta^t_j}, ~\forall t\in [T], j\in [n], k\in\{0,\ldots,K\}
\end{align*}
Here we use the same notation \(K\) for all \(t \in \{1, \ldots, T\}\) and \(j \in [n]\), to simplify the notations. However, it's important to note that \(K\) can vary across different stages and products.

As discussed in the previous section, solving [\textbf{MILP}$^t_{i}$] can provide solutions within an \(\mathcal{O}(1/K)\) neighborhood of the true solution. By repeatedly solving [\textbf{MILP}$^t_{i}$] for \(t = T, T-1, \ldots, 1\) and \(j \in [n]\), we can approximately solve [\textbf{DP}]. The main advantage of this discretization and approximation approach, compared to prior works, is that it allows handling any linear constraints on both prices and purchasing probabilities (while prior works only handle constraints on the purchasing probabilities). Although this approach requires repeatedly solving MILPs, which are typically more expensive than solving convex programs when there are no constraints on the price variables, state-of-the-art MILP solvers (such as CPLEX and GUROBI) can solve large MILPs in seconds with support from advanced machines. In our later experiments, we demonstrate that our approach not only provides better solutions than other baselines but can also solve large instances in a reasonable computing time.

\section{Baselines}
We provide a more detailed description of the baseline methods considered in our experiments.
\subsubsection{Static Pricing:}
For solving [\textbf{SP}], we will compare our method described in Section \ref{sec:SP} with some standard baselines in the literature. Specifically, we will consider the following methods for solving the  static pricing problem:
\begin{itemize}
    \item \texttt{SP-DMIP:} This is our method, standing for Discretization and MILP. Details of this approach can be found in Section \ref{sec:SP}.    
    \item \texttt{SP-SLSQP:} In this approach, we use the \textit{Sequential Least Squares Programming} method \citep{kraft1988software} to directly solve the \textbf{[SP]}. This is a standard and popular method to solve nonlinear optimization problems.
    \item \texttt{SP-Trans:} As discussed earlier, prior pricing methods typically transform the pricing problem into a convex optimization program with constraints on the purchasing probabilities. We adapted such an approach to compare with our method. Specifically, according to \citep{dong2009dynamic,zhang2013assessing}, when the price sensitivity parameters are the same over products, one can see that \(P_j/P_0 = \exp((a_j-r_j)/b_j)\). Then the price of each product \(r_j\) can be written as a function of purchase probabilities \(\bP\) where 
    \begin{equation}
    \label{prob transfer}
        r_j(\bP) = a_j - b_j\ln P_j + b_j \ln P_0 \quad \forall j \in [n]   
    \end{equation}
    Thus, instead of using the vector \(\br\) as decision variables for \textbf{[SP]}, we can perform the above change of variables and use the vector \(\bP\) as decision variables. The reformulation of \textbf{[SP]} is given as follows:
    \begin{align}
    [\textbf{SP}_p] ~~~~~  \max \limits_{\bP\geq \textbf{0}}  &\left\{ d \sum_{j=1}^n P_j(a_j - b_j\ln P_j + b_j \ln P_0) \right\} \nonumber \\
    \mbox{subject to} &\quad d\bA \bP \leq \bc   \label{SP_p_capacity}  \\
    & \quad \sum_{j=1}^n P_{j} + P_0  = 1   \label{SP_p_prob} \\
    & \quad  \bB\br(\bP) \leq \bd  \label{SP_p_addtion2} \\
    & \quad  \br_j(\bP)\in [l_j,u_j] ~~\forall j\in [n] \label{SP_p_addtion3}
    \end{align}
    It is well-known that, without the price constraints in \eqref{SP_p_addtion2} and \eqref{SP_p_addtion3}, and the price-sensitivity parameters are the same  across product, i.e., $b_j=b_{j'}, ~\forall j,j'\in [n]$, \textbf{[SP$_p$]} is a concave maximization problem, which can be efficiently solved to optimality by a standard convex solver. To adapt this to handle pricing constraints, we first remove Constraints \eqref{SP_p_addtion2} and \eqref{SP_p_addtion3} from \textbf{[SP$_p$]} and solve the resulting convex optimization problem to obtain a purchasing vector \(\bP^*\) and corresponding pricing solution \(\br^*\) (using equation (\ref{prob transfer})). Finally, to get a feasible pricing solution, we project \(\br^*\) into the feasible set by solving the following problem:
    \begin{align}
    ~~~~ \min \limits_{\br} & \quad \lVert \br-\br^* \rVert  \\
    \mbox{subject to} &  \quad \texttt{Constraints } (\ref{SP_p_addtion2} - \ref{SP_p_addtion3}) 
    \end{align}
    
    \item \texttt{Gurobi:} Recently, GUROBI \citep{GUROBI} has become a state-of-the-art solver for mixed-integer nonlinear programs, capable of handling nonlinear components by piece-wise linear approximation (PWLA) approaches. While GUROBI's PWLA cannot directly solve [\textbf{SP}] due to the fractional structure, we can apply it to solve each step of the bisection described in Section \ref{sec:SP}:
    \begin{align}
    ~~~~~ \max_{\br \in \cR} & \quad \sum_{j\in[n]} f_j(r_j) - \delta \left(\sum_{j\in [n]} g_j(r_j) + 1\right)
    \end{align}
    where \(f_j(r_j) = r_j e^{(a_j-r_j)/b_j}\) and \(g_j(r_j) = e^{(a_j-r_j)/b_j}\). To handle these two nonlinear univariate functions, we approximate \(f\) and \(g\) with piece-wise linear functions. In the GUROBI implementation, we add new variables involved in the piece-wise linear function, define the breakpoints of \(f_j\) and \(g_j\) along with the price \(r_j\), and transform the model into a MILP. To have a fair comparison, we use the same number of breakpoints \(K\) as in our approach. So, in short, the \texttt{Gurobi} approach involves breaking [\textbf{SP}] into a sequence of subproblems, and solving each subproblem by PWLA implemented in the solver with the same number of breakpoints.
\end{itemize}


\subsubsection{{Dynamic Pricing:}}
We employ the dynamic programming decomposition approach described in Section \ref{section DP}. After solving all the subproblem [\textbf{DP}$_i$], we can collect all the approximate value function  $\{ v_{t,i}(\cdot)\}, {\forall {t\in[T],i\in [m]}}$. The value function $v_t(\bx)$ can then be approximated by
\begin{equation}
\label{dp_decomposition of value function1}
v_t(\bx) \approx \sum_{i=1}^m v_{t,i}(x_i)
\end{equation}
Then, a  pricing solution can be computed accordingly as:
\begin{equation}
\label{dp policy}
\br^*_t(\bx) = \mathop{\arg \max} \limits_{\br_t \in \cR_t(\bx)}~~~\left\{\sum_{j=1}^n \lambda P_j(\br_t)[r_{t,j} - \sum_{i=1}^m \triangle_j v_{t+1,i}(x_i)] \right\}
\end{equation}
where
\begin{equation}
\label{delta value function1}
\triangle_j v_t(\bx) = v_t(\bx) - v_t(\bx-\bA^j) \approx \sum_{i=1}^m \triangle_j v_{t,i}(x_i),~\forall t\in \{1,\ldots,T+1\},~ j\in [n]
\end{equation}
We compare our approach, which solves [\textbf{DP}$_i$] using our discretization and MILP approximation approach, with several baselines, similar to the static pricing case. These methods are described as follows:

\begin{itemize}
    \item \texttt{DP-DMIP}: This is our method, solving \textbf{[DP$_i$]} using the discretization and MILP approximation described in Section \ref{sec:DP}.
    
    \item \texttt{DP-SLSQP}: Similar to static pricing, this approach relies on the use of the \textit{Sequential Least Squares Programming }method to directly solve \textbf{[DP$_i$]}.
    
    \item \texttt{DP-Trans}: We use the same method as \texttt{SP-Trans} for solving \textbf{[DP$_i$]}. Specifically, by considering the probability vector \(\bP\) as decision variables, we first remove the constraints on the prices and solve the concave maximization problem at each stage:
    
    \begin{equation}
    \label{DP_p}
    [\textbf{{DP}}_{ip}] \quad \quad \max \limits_{\bP_t \geq \textbf{0}}  \lambda \sum_{j=1}^n P_{t,j} \left( a_j - b_j\ln P_{t,j} + b_j \ln P_{t,0} - \sum_{k \neq i} a_{kj}\pi_k - \triangle_j v_{t+1,i}(x_i) \right)
    \end{equation}
    
    The obtained probability solution \(\bP^*_t\) can be used to retrieve a pricing solution \(\br^*_t\), which is then projected to the feasible set \(\cR_t(\bx)\) to get a feasible pricing solution.
\end{itemize}




\subsection{Choosing the Number of Breakpoints $K$:}

\begin{figure}[htb]
    \begin{subfigure}{0.5\textwidth} 
        \centering
        \includegraphics[width=\linewidth]{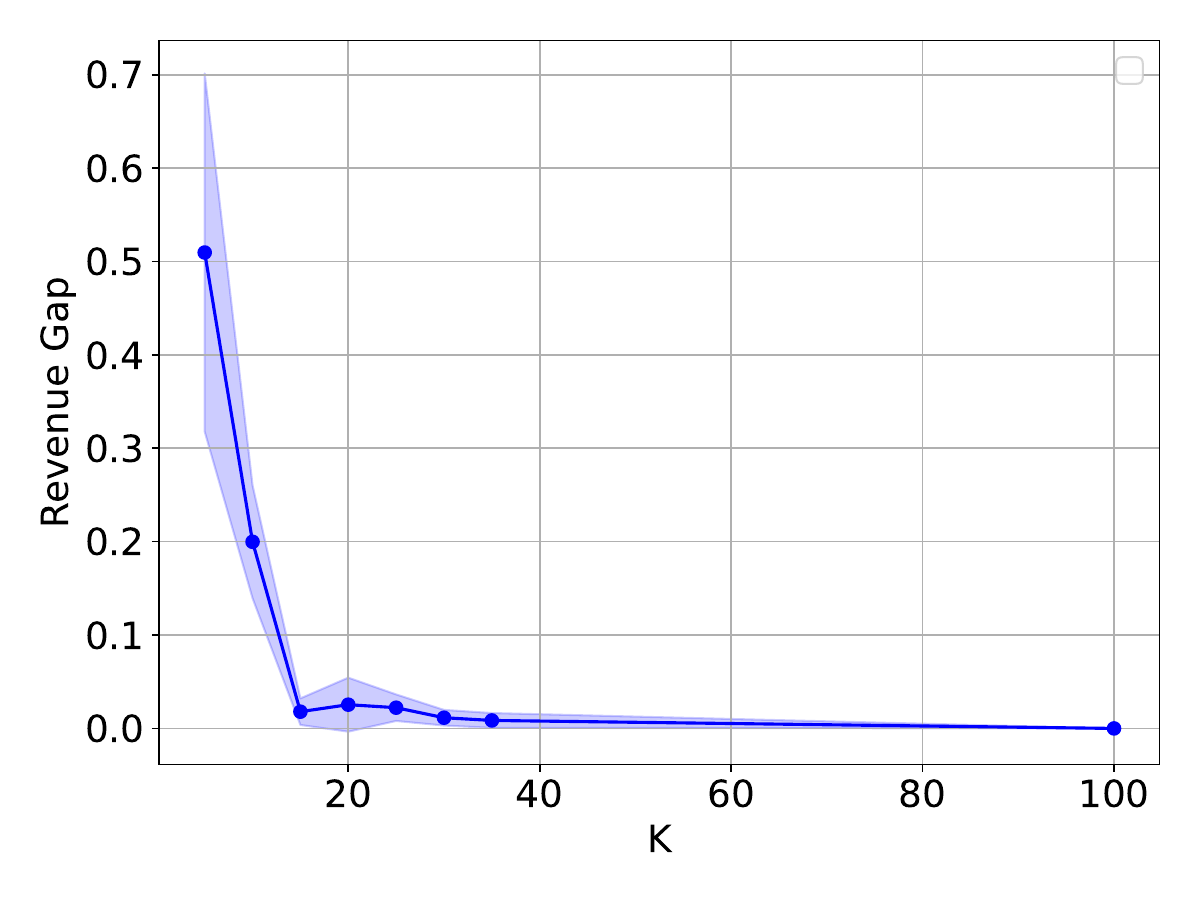} 
        \caption{Revenue Gap \%}
        \label{fig:Revenue Gap}
    \end{subfigure}%
    \begin{subfigure}{0.5\textwidth} 
        \centering
        \includegraphics[width=\linewidth]{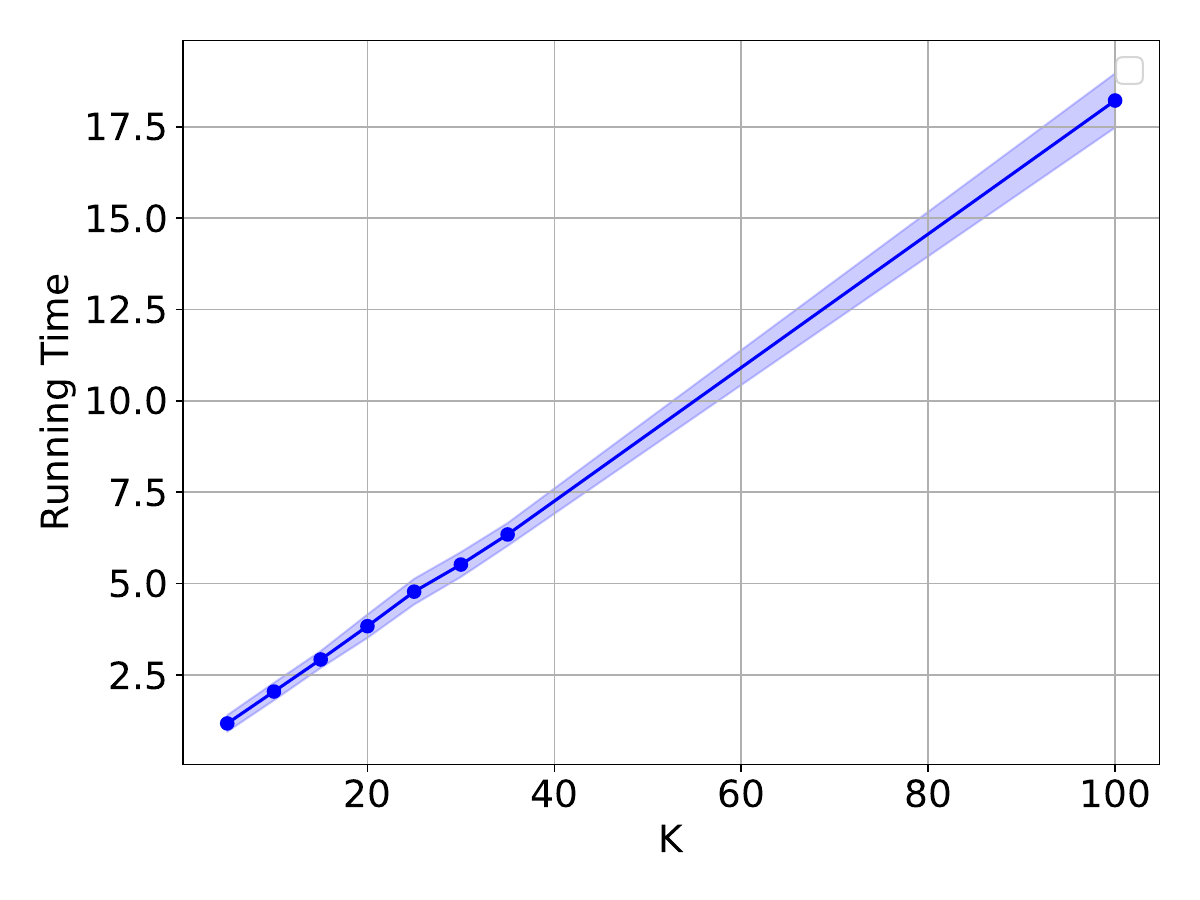} 
        \caption{Running Time}
        \label{fig:Running Time}
    \end{subfigure}
    \caption{For m=16, n=80, T=200, the mean and standard deviation of revenue gap/running time of 5 different random seeds }
    \label{fig:rev and time}
\end{figure}

 A crucial aspect of our piecewise linear approximation approach is the number of breakpoints \( K \) used to discretize the price interval \([l_j, u_j]\) for each product \( j \). As shown in Section \ref{sec:SP}, the approximation error is inversely proportional to \( K \), whereas the sizes of the MILP formulations increases with \( K \). Therefore, selecting an optimal \( K \) is essential to balance between the accuracy of the approximation and the computational size of the model. To this end, we numerically assess the approximation errors as a function of \( K \). We randomly generate five instances and solve them with different values of \( K \). For each \( K \), we calculate the revenue gap as the percentage difference between the \textit{near-optimal revenue} (achieved when \( K \) is sufficiently large, i.e., \( K=100 \)) and the revenue obtained with that \( K \):
\[
\frac{\texttt{[Near-optimal revenue]} - [\texttt{Revenue given by } K]}{\texttt{[Near-optimal revenue]}} \times 100\%
\]
Figure \ref{fig:rev and time} shows the changes in the revenue gap and running time as the parameter \( K \) increases. The graph indicates that when \( K \) exceeds 15, the revenue gap approaches zero, suggesting that \( K=15 \) is sufficient to achieve near-exact objective function values. Additionally, the running time linearly increases with higher values of \( K \). Based on these findings, we have opted to set \( K \) at 15 for our later experiments.

\end{APPENDICES}

\end{document}